\documentclass[11pt,
]{amsart}
\usepackage{
amssymb, graphicx
}
\usepackage{mathrsfs} 
\usepackage{amsthm}
\usepackage{graphicx,color}
\usepackage[dvipsnames]{xcolor}
\usepackage{amsmath}
\allowdisplaybreaks[4]
\newtheorem{theorem}{Theorem}
\newtheorem{lemma}{Lemma}
\newtheorem{proposition}{Proposition}

\newtheorem{corollary}{Corollary}
\newtheorem{remark}{Remark}
\theoremstyle{remark}

\date{\today}

\title[Determination of the density]{Determination of the density in the linear elastic wave equation}

  \author[J. Zhai]{Jian Zhai}
\address{School of Mathematical Sciences,
  Fudan University, Shanghai 200433, China
  (\tt{jianzhai@fudan.edu.cn}).}
    \thanks{J. Zhai is supported by National Key Research and Development Programs of China (No. 2023YFA1009103), NSFC(No. 12471396), Science and Technology Commission of Shanghai Municipality (23JC1400501)}

\begin{document}
\begin{abstract}
We study the inverse boundary value problem for the linear elastic wave equation in three-dimensional isotropic medium. We show that both the Lam\'e parameters and the density can be uniquely recovered from the boundary measurements under the strictly convex foliation condition.
\end{abstract}
\maketitle

\section{Introduction}
The wave equation for the displacement vector $u=(u_1,u_2,u_3)$ in a three-dimensional isotropic elastic body $\Omega$ reads
\begin{equation}\label{elastic_eq}
(Pu)_i=\rho\partial^2_t u_i-\sum_{jk\ell}\partial_j(c_{ijk\ell}\partial_\ell u_k)=0,
\end{equation}
where $c_{ijk\ell}=\lambda\delta_{ij}\delta_{k\ell}+\mu(\delta_{ik}\delta_{j\ell}+\delta_{i\ell}\delta_{jk})$. Here we assume that the Lam\'e parameters $\lambda,\mu$ and the density $\rho$ all belong to $C^\infty(\overline{\Omega})$ and
\begin{equation}\label{positivitycondition}
\mu>0,\quad 3\lambda+2\mu>0,\quad \rho>0,\quad\text{on }\overline{\Omega}.
\end{equation}
 We will consider the initial boundary value problem
\begin{equation}
\begin{cases}
Pu=0,\quad &\text{in } (0,T)\times \Omega,\\
u=f,\quad &\text{on }(0,T)\times\partial\Omega,\\
u=\partial_tu=0,\quad &\text{on }\{t=0\}\times\Omega,
\end{cases}
\end{equation}
where $\Omega\subset\mathbb{R}^3$ is an open domain, and the boundary $\partial\Omega$ is smooth. We will consider the inverse boundary value problem of recovering $\lambda,\mu,\rho$ from the Dirichlet-to-Neumann map 
\[
\Lambda_{\lambda,\mu,\rho}:f\mapsto \sum_{jk\ell}\nu_jc_{ijk\ell}\partial_\ell u_k\vert_{(0,T)\times\partial\Omega},
\]
where $\nu$ is the unit outer normal to $\partial\Omega$. The main application of this problem is seismic inversion \cite{baeten1989theoretical,shi2020numerical}, that is to recover the subsurface geologic structure from seismic data.\\

The principal symbol of $P$, as a pseudodifferential operator, is
\[
p_{ik}(t,x,\tau,\xi)=\rho(x) \tau^2\delta_{ik}-c_{ijk\ell}(x)\xi_j\xi_\ell.
\]
The characteristic set $\{(t,x,\tau,\xi)\vert \det p(t,x,\tau,\xi)=0\}$ of $P$ consists of two separate sets $\Sigma_P=\{(t,x,\tau,\xi)\vert \tau^2-c_P^2(x)|\xi|^2=0\}$ and $\Sigma_S=\{(t,x,\tau,\xi)\vert \tau^2-c_S^2(x)|\xi|^2=0\}$, where
\[
c_P=\sqrt{\frac{\lambda+2\mu}{\rho}},\quad\quad c_S=\sqrt{\frac{\mu}{\rho}}
\]
are the wavespeeds for $P$- and $S$- waves. Clearly \eqref{positivitycondition} implies $c_P>c_S$. The $P/S$-wave paths are then modeled by integral curves (zero-bicharacteristics) $(t(s),x(s),\tau(s),\xi(s))$ of the Hamiltonian vector field associated with the Hamiltonian $\tau^2-c^2_{P/S}(x)|\xi|^2$. The bicharacteristics $(x(s),\xi(s))$ is just the lift-ups of geodesics with respect to the Riemannian metric $g_{P/S}:=c_{P/S}^{-2}g_E$ to the phase space, where $g_E$ is the Euclidean metric.

The corresponding problem for the scalar wave equation $(\partial^2_t-\Delta_g)u=0$ is well studied, where $g$ is a Riemannian metric. It is well known that one can recover $g$ up to an isometry fixing $\partial\Omega$ using the boundary control method \cite{belishev1992reconstruction}. For a different approach, under certain geometric conditions, the problem can also be reduced to some geometric inverse problems, and then uniqueness and stability results can be derived \cite{bao2014sensitivity,stefanov2005stable,stefanov2016stable}.

For the elastic wave equations, the inverse problem is less understood. Rachele \cite{rachele2000boundary} proved that the jets of $\lambda,\mu,\rho$ at the boundary can be fully recovered. For interior determination, the boundary control method has not been successfully applied to the full Lam\'e equation \eqref{elastic_eq}. For a Lam\'e-type equation having the same principal part that admits a full $P/S$-decoupling, the boundary control was used to recover the wavespeeds \cite{belishev2006dynamical,belishev2013reachable,fomenko2021determination}. Via reduction to geometric inverse problems, the uniqueness of the two wavespeeds $c_P,c_S$ (or equivalently $\frac{\lambda}{\rho},\frac{\mu}{\rho}$) is known (1) assuming $(\Omega,g_{P/S})$ is simple \cite{rachele2000inverse}; or (2) assuming that $(\Omega,g_{P/S})$ satisfies the strictly convex foliation condition \cite{stefanov2017local}. Recall that a compact Riemannian manifold $(\Omega,g)$ with strictly convex boundary $\partial\Omega$ satisfies the strictly convex foliation condition if it admits a smooth strictly convex function (see \cite{uhlmann2016inverse,paternain2019geodesic} for more discussions on this condition). We note here that a function $h$ is called strictly convex if each level set $h=c$ is strictly convex when viewed from $h>c$. The proof is based on recovering the lens relations related to the two wavespeeds and then applying the known rigidity results \cite{croke1991rigidity,stefanov2016boundary}. With $c_P$ and $c_S$ already determined, under the above two geometric assumptions, \cite{rachele2003uniqueness,bhattacharyya2018local} study the recovery of the density $\rho$ using the subprincipal term of the $P$-wave amplitude  and the inversion of the longitudinal geodesic ray transform of a second order tensor.  But the results need to essentially assume a technical assumption $\lambda\neq2\mu$. These results can be extended to uniqueness of piecewise smooth parameters \cite{caday2021recovery,bhattacharyya2022recovery}. We also refer to \cite{chen2024stability} for a uniqueness result for the density under different assumptions.  For a nonlinear elastic wave equation, the density can be recovered utilizing the nonlinearity \cite{uhlmann2024determination}.\\

In this work we study the recovery of the density $\rho$ using $S$-wave amplitude with known $c_P$ and $c_S$. With $c_{P/S}$ already known, we use the notation $\Lambda_\rho=\Lambda_{\lambda,\mu,\rho}$ for simplicity. The readers will see that the problem reduces to the inversion of the transverse geodesic ray transform of a second order tensor. To be more specific, we consider a Riemannian manifold $(\Omega,g)$ with boundary $\partial\Omega$. The transverse geodesic ray transform of a second order tensor $f$ is
\[
T_2f(\gamma,\eta):=\int_\gamma \langle f(\gamma(s)),\eta\otimes\eta\rangle_g\mathrm{d}s=\int_\gamma g^{ik}g^{j\ell}f_{ij}(\gamma(s))\eta_k(s)\eta_\ell(s)\mathrm{d}s,
\]
where $\gamma$ runs over any geodesic in $(\Omega,g)$ connecting two boundary points, and $\eta$ is an arbitrary parallel vector field along $\gamma$ and orthogonal to $\dot{\gamma}$. In contrast, the longitudinal geodesic ray transform, encountered in \cite{rachele2003uniqueness,bhattacharyya2018local}, of $f$ is
\[
I_2f(\gamma):=\int_\gamma \langle f(\gamma(s)),\dot{\gamma}\otimes\dot{\gamma}\rangle_g\mathrm{d}s.
\]
         
The injectivity of $I_2$ has been studied quite extensively (cf. \cite{Shara,stefanov2004stability,dairbekov2006integral,paternain2015invariant,stefanov2018inverting}). For the transverse ray transform, Sharafutdinov \cite{Shara} proved the injectivity of $T_2$ assuming that $(\Omega,g)$ is simple and the curvature is close to zero.
 Recently, the author and collaborator \cite{uhlmann2024invertibility1,uhlmann2024invertibility2} proved the injectivity of $T_2$ under the strictly convex foliation condition.
\begin{proposition}\label{thmglobal20}{\textnormal{(\cite{uhlmann2024invertibility2})}}
Assume $\partial \Omega$ is strictly convex and $(\Omega,g)$ admits a smooth strictly convex function. Assume also that the second order tensor $f$ is symmetric. If
\[
T_2f(\gamma,\eta)=0
\]
for any geodesic $\gamma$ in $\Omega$ with endpoints on $\partial \Omega$, and $\eta$ a parallel vector field orthogonal to $\gamma$. Then $f=0$.
\end{proposition}
\noindent As pointed out in \cite{stefanov2016boundary} the foliation condition is a natural generalization of the Herglotz \cite{herglotz1905elastizitaet} and the Wieckert-Zoeppritz \cite{wiechert1907erdbebenwellen} conditions, which are typically assumed for the Earth. 

We will prove that from the Dirichlet-to-Neumann map one can recover $T_2N$ for some tensor $N$ depending on $\rho$, where the Riemannian metric $g=g_S$ is induced by the $S$-wavespeed. This is done by the construction and analysis of geometric optics solutions, and the explicit form of $N$ is given in \eqref{formofN}. Under the foliation condition, the invertibility of $T_2$ (cf. Proposition \ref{thmglobal20}) guarantees that we can recover the symmetrization of $N$, $\mathrm{sym}(N)$. Finally we will then recover $\rho$ from $\mathrm{sym}(N)$ by solving a fourth order elliptic partial differential equation.\\

The main result of this paper is as follows.

\begin{theorem}\label{maintheorem}
Let $c_P$ and $c_S$ be already known. Assume $\partial\Omega$ is strictly convex with respect to $g_{S}$, and $(\Omega,g_{S})$ admits a strictly convex function.
If $T$ is greater than the length of all geodesics in $(\Omega,g_{S})$, then $\Lambda_{\rho}$ determines $\rho$ in $\overline{\Omega}$ uniquely.
\end{theorem}

\begin{remark}
\textnormal{Compared with the results \cite{rachele2003uniqueness,bhattacharyya2018local}, we have removed the artificial assumption $\lambda\neq2\mu$.}
\end{remark}

Under the strictly convex foliation conditions for both $g_P$ and $g_S$,  the uniqueness results for $\frac{\lambda}{\rho}$ and $\frac{\mu}{\rho}$ have been obtained in \cite{stefanov2017local}. Actually, only the singular part of the Schwartz kernel of the Dirichlet-to-Neumann map is needed.
\begin{proposition}\label{uniquenessspeeds}{\textnormal{(\cite{stefanov2017local})}}
Assume $\partial\Omega$ is strictly convex with respect to $g_{P/S}$, and $(\Omega,g_{P/S})$ admits a strictly convex function. If $T$ is greater than the length of all geodesics in $(\Omega,g_{S})$, then the knowledge of the Schwartz kernel of $\Lambda_{\lambda,\mu,\rho}$, up to a smooth function, determines $c_{P/S}$ in $\overline{\Omega}$ uniquely.
\end{proposition}

With the above proposition, we have the following corollary of Theorem \ref{maintheorem}.
\begin{corollary}
 Assume $\partial\Omega$ is strictly convex with respect to $g_{P}$ and $g_{S}$, and $(\Omega,g_{P/S})$ admits a strictly convex function.
If $T$ is greater than the length of all geodesics in $(\Omega,g_{P})$ and $(\Omega,g_{S})$, then $\Lambda_{\lambda,\mu,\rho}$ determines $\lambda,\mu,\rho$ in $\overline{\Omega}$ uniquely.
\end{corollary}

\begin{remark} 
\textnormal{One can also apply the uniqueness result of $T_2$ in \cite{Shara} to obtain the uniqueness of $\rho$ under the assumption that the Riemannian manifold $(\Omega,g_S)$ is simple and its curvature is close to zero.}
\end{remark}

\begin{remark}
\textnormal{As in \cite{stefanov2017local,bhattacharyya2018local}, one can also prove a local uniqueness result, namely, for recovering the parameter near a convex part of the boundary from measurements taken only on that part.}
\end{remark}
 
 The rest of the paper is organized as follows. In Section \ref{GOsolutions}, we outline the procedure for the construction of geometric optics solutions. In Section \ref{Pwave}, we use $P$-wave geometric optics solutions to reproduce the results in \cite{rachele2003uniqueness,bhattacharyya2018local}. In Section \ref{Swave}, we use $S$-wave geometric optics solutions to obtain the main result. The explicit construction of the geometric optics solutions used in Section \ref{Pwave} and Section \ref{Swave} will be provided in the appendix.
 
\section{Geometric optics solutions}\label{GOsolutions}
The proof relies on the construction of geometric optics solutions expressed using Fourier integral operators \cite{hormander1971fourier,duistermaat1994fourier}. Assume $V$ is an open subset of $\partial\Omega$.
We will construct such solutions to the elastic wave equation \eqref{elastic_eq} of the form
\begin{equation}\label{solutiondistribution}
u(t,x)=(2\pi)^{-1}\int_\mathbb{R} e^{\mathrm{i}\varrho(t-\varphi_{P/S}(x))}\mathbf{a}(t,x,\varrho)\mathrm{d}\varrho,
\end{equation}
where $\mathbf{a}\in S^m$ is a vector. The wavefront set of the above distribution is contained in
\begin{equation}\label{wavefrontset}
\{(t,x,\varrho,-\varrho\nabla_x\varphi_{P/S});t=\varphi_{P/S}(x),\varrho\neq 0\}.
\end{equation}
The phase functions $\varphi_{P/S}(x)$ solve the Eikonal equation
\begin{equation}\label{eikonal0}
c^2_{P/S}|\nabla_x\varphi_{P/S}|^2=1.
\end{equation}
The solution to the Eikonal equation exists only locally in general.
Then the wave paths are given by
\[
(\tau, \xi)=\varrho(1,-\nabla_x\varphi_{P/S}).
\]

Fix $(t_0,x_0,\tau_0,\xi_0)$, $x_0\in\partial\Omega$, in the set \eqref{wavefrontset}. Take $V\subset \partial\Omega$ a neighborhood of $x_0$. We can microlocalize $u$ such that it is defined near the ($P$- or $S$-wave) null-geodesic $\vartheta_0$ issued from $(t_0,x_0)\in (0,T)\times V$ in the direction $(\tau_0,\xi_0)$ but in some neighborhood of $(t_0,x_0)$. To extend it globally, we assume that it is valid until $\vartheta_0$ hits another hypersurface $(0,T)\times V_1$, and a bit beyond it. We then consider a new similar problem with $V$ replaced by $V_1$. By a compactness argument, we can cover the whole null-geodesic.

Assume that $\mathbf{a}$ admits an asymptotic expansion as $\varrho\rightarrow\infty$
\[
\mathbf{a}(t,x,\varrho)\sim\sum_{j=0}^\infty\varrho^{-j} \mathbf{a}_j(t,x).
\]
We will carry out the construction of $\mathbf{a}$ in the so-called ray coordinates (associated with $g_\bullet=c^{-2}_\bullet\mathrm{d}s^2$) for $\mathbb{R}^3$, i.e., curvilinear coordinates $(x^1,x^2,x^3)$ such that $x^3(x)=\tau=\varphi_\bullet(x)$ and the coordinate surface $x^3=x^3_0$ are orthogonal to the coordinate lines $x^1=x_0^1$, $x^2=x^2_0$ that are the geodesics of metric $g_\bullet$. Here $\bullet=P$ or $S$, and we occasionally omit it.  Under these ray coordinates the Euclidean metric has the form
\[
g_E=\mathrm{d}s^2=g_{\alpha\beta}\mathrm{d}x^\alpha\mathrm{d}x^\beta+c^2_\bullet\mathrm{d}\tau^2.
\]
From now on the Greek indices range over the values $1,2$.
Throughout the paper, we treat $g_E,g_P,g_S$ as covariant second order tensors. Under ray coordinates, the Euclidean metric is expressed as $(g_{ij})$ with $g_{3j}=c_\bullet^2\delta_{3j}$, $i,j=1,2,3$. The matrix $(g^{ij})$ is the inverse matrix of $(g_{ij})$. Then obviously, under the ray coordinates the metric $g_\bullet=c_\bullet^{-2}g_E$ is expressed as  $(c_\bullet^{-2} g_{ij})$.
The Christoffel symbols of the Euclidean metric are
\[
\begin{split}
&\Gamma^3_{\alpha\beta}=-\frac{1}{2}c^{-2}\frac{\partial g_{\alpha\beta}}{\partial\tau},\quad \Gamma^\beta_{\alpha 3}=\frac{1}{2}g^{\beta\delta}\frac{\partial g_{\alpha\delta}}{\partial \tau},\quad \Gamma^3_{\alpha 3}=c^{-1}\frac{\partial c}{\partial x^\alpha},\\
&\Gamma^\alpha_{33}=-cg^{\alpha\beta}\frac{\partial c}{\partial x^\beta},\quad \Gamma^3_{33}=c^{-1}\frac{\partial c}{\partial \tau}.
\end{split}
\]
Note that the equation \eqref{elastic_eq} can be written as
\[
\rho \partial^2_tu-\nabla\cdot\sigma(u)=0,
\]
where
\[
\varepsilon(u)=\frac{1}{2}(\nabla u+(\nabla u)^T)
\]
is the strain tensor and
\[
\sigma(u)=\lambda\mathrm{trace}(\varepsilon(u))e+2\mu\varepsilon(u)
\]
is the stress tensor. Here $e$ is the Euclidean metric tensor.

For ease of notation we write $\varphi=\varphi_P$ or $\varphi_P$. Denote
\[
u_\varrho=\mathbf{a}e^{\mathrm{i}(t-\varrho\varphi)},
\]
and write $\mathbf{a}=(a_{1},a_{2},a_{3})$ under the ray coordinates.
The strain tensor takes the form
\[
\varepsilon_{k\ell}(u_\varrho)=\frac{1}{2}\mathrm{i}\varrho(a_k\varphi_{;\ell}+a_\ell\varphi_{;k})e^{\mathrm{i}\varrho\varphi}+\frac{1}{2}(a_{k;\ell}+a_{\ell;k})e^{\mathrm{i}\varrho\varphi}.
\]
The stress tensor is then
\[
\begin{split}
\sigma_{ij}(u_\varrho)=&\lambda e^{k\ell}\varepsilon_{k\ell}e_{ij}+2\mu\varepsilon_{ij}\\
=&\lambda g^{k\ell}\varepsilon_{k\ell} g_{ij}+2\mu \varepsilon_{ij}\\
=&\mathrm{i}\varrho(\lambda a_k\varphi_{;\ell}g^{k\ell}g_{ij}+\mu a_i\varphi_{;j}+\mu a_j\varphi_{;i})e^{\mathrm{i}\varrho\varphi}+(\lambda a_{k;\ell}g^{k\ell}g_{ij}+\mu a_{i;j}+\mu a_{j;i})e^{\mathrm{i}\varrho\varphi}.
\end{split}
\]
We calculate
\begin{align*}
\sigma_{ij;m}(u_\varrho)=&\partial_m\sigma_{ij}-\Gamma^n_{im}\sigma_{nj}-\Gamma^n_{jm}\sigma_{ni}\\
=&-\varrho^2(\lambda a_k\varphi_{;\ell}g^{k\ell}g_{ij}+\mu a_i\varphi_{;j}+\mu a_j\varphi_{;i})\varphi_{;m}e^{\mathrm{i}\varrho\varphi}\\
&+\mathrm{i}\varrho\partial_m(\lambda a_k\varphi_{;\ell}g^{k\ell}g_{ij}+\mu a_i\varphi_{;j}+\mu a_j\varphi_{;i})e^{\mathrm{i}\varrho\varphi}\\
&+\mathrm{i}\varrho(\lambda a_{k;\ell}g^{k\ell}g_{ij}+\mu a_{i;j}+\mu a_{j;i})\varphi_{;m}e^{\mathrm{i}\varrho\varphi}\\
&-\mathrm{i}\varrho\Gamma^n_{im}(\lambda a_k\varphi_{;\ell}g^{k\ell}g_{nj}+\mu a_n\varphi_{;j}+\mu a_j\varphi_{;n})e^{\mathrm{i}\varrho\varphi}\\
&-\mathrm{i}\varrho\Gamma^n_{jm}(\lambda a_k\varphi_{;\ell}g^{k\ell}g_{ni}+\mu a_n\varphi_{;i}+\mu a_i\varphi_{;n})e^{\mathrm{i}\varrho\varphi}\\
&+\partial_m(\lambda a_{k;\ell}g^{k\ell}g_{ij}+\mu a_{i;j}+\mu a_{j;i})e^{\mathrm{i}\varrho\varphi}\\
&-\Gamma_{im}^n(\lambda a_{k;\ell}g^{k\ell}g_{nj}+\mu a_{n;j}+\mu a_{j;n})e^{\mathrm{i}\varrho\varphi}\\
&-\Gamma_{jm}^n(\lambda a_{k;\ell}g^{k\ell}g_{ni}+\mu a_{n;i}+\mu a_{i;n})e^{\mathrm{i}\varrho\varphi}.
\end{align*}

Therefore
\begin{align*}
\sigma_{ij;}^{~\,~j}(u_\varrho)=&\sigma_{ij;m}g^{jm}\\
=&-\varrho^2(\lambda a_k\varphi_{;\ell}g^{k\ell}\varphi_{;i}+\mu a_j\varphi_{;m}g^{jm}\varphi_{;i}+\mu\varphi_{;j}\varphi_{;m}g^{jm}a_i)e^{\mathrm{i}\varrho\varphi}\\
&+\mathrm{i}\varrho\left(\partial_i(\lambda a_k\varphi_{;\ell}g^{k\ell})+\lambda a_k\varphi_{;\ell}g^{k\ell}\partial_mg_{ij} g^{jm}+\partial_m(\mu a_i\varphi_{;j}+\mu a_j\varphi_{;i})g^{jm}\right)e^{\mathrm{i}\varrho\varphi}\\
&+\mathrm{i}\varrho\left(\lambda a_{k;\ell}g^{k\ell}\varphi_{;i}+\mu(a_{i;j}+a_{j;i})\varphi_{;m}g^{jm}\right)e^{\mathrm{i}\varrho\varphi}\\
&-\mathrm{i}\varrho\Gamma^n_{im}(\lambda a_k\varphi_{;\ell}g^{k\ell}g_{nj}+\mu a_n\varphi_{;j}+\mu a_j\varphi_{;n})g^{jm}e^{\mathrm{i}\varrho\varphi}\\
&-\mathrm{i}\varrho\Gamma^n_{jm}(\lambda a_k\varphi_{;\ell}g^{k\ell}g_{ni}+\mu a_n\varphi_{;i}+\mu a_i\varphi_{;n})g^{jm}e^{\mathrm{i}\varrho\varphi}\\
&+\partial_m(\lambda a_{k;\ell}g^{k\ell}g_{ij}+\mu a_{i;j}+\mu a_{j;i})g^{jm}e^{\mathrm{i}\varrho\varphi}\\
&-\Gamma_{im}^n(\lambda a_{k;\ell}g^{k\ell}g_{nj}+\mu a_{n;j}+\mu a_{j;n})g^{jm}e^{\mathrm{i}\varrho\varphi}\\
&-\Gamma_{jm}^n(\lambda a_{k;\ell}g^{k\ell}g_{ni}+\mu a_{n;i}+\mu a_{i;n})g^{jm}e^{\mathrm{i}\varrho\varphi}.
\end{align*}
Here $\langle\cdot,\cdot\rangle$ denotes the Euclidean inner product and $|\cdot|$ denote the Euclidean norm.
Note that
\[
\rho\partial^2_tu_\varrho=-\varrho^2\rho\mathbf{a}e^{\mathrm{i}\varrho(t-\varphi)}+2\mathrm{i}\varrho \rho(\partial_t\mathbf{a})e^{\mathrm{i}\varrho(t-\varphi)}+\rho\partial^2_t\mathbf{a}e^{\mathrm{i}\varrho(t-\varphi)}.
\]

In order to construct microlocal solutions to \eqref{elastic_eq}, we need to take $u_\varrho$ such that
\begin{equation}\label{elastic_eq_cova}
\rho\partial^2_t (u_\varrho)_i-\sigma_{ij;}^{~\,~j}(u_\varrho)\sim 0.
\end{equation}
One can use the above equation to derive a sequence of equations for $\varphi$ and $\mathbf{a}_j$ by equating the coefficients of same orders in $\varrho$. We refer to \cite{rachele2000inverse,rachele2003uniqueness} for a construction of microlocal solutions under Cartesian corrdinates.

Equating the coefficients of order $\varrho^2$ in the equation \eqref{elastic_eq_cova}, we get
\begin{equation}\label{first_eq}
\rho\mathbf{a}_0=(\lambda+\mu)\langle\mathbf{a}_0,\nabla\varphi\rangle\nabla\varphi+\mu|\nabla \varphi|^2\mathbf{a}_0.
\end{equation}
For $P$ waves, the wave vector $\nabla\varphi$ and the principal amplitude/polarization $\mathbf{a}_0$ are parallel. Actually, if we take $\mathbf{a}_0=a\nabla\varphi$ for some scalar function $a$, then the equation \eqref{first_eq} becomes
\[
\rho\mathbf{a}_0=(\lambda+2\mu)|\nabla\varphi|^2\mathbf{a}_0,
\]
which is automatically satisfied by choosing $\varphi=\varphi_P$ satisfying $c_P^2|\nabla\varphi_P|^2=1$ (cf. \eqref{eikonal0}).
For $S$-waves, we take $\langle\mathbf{a}_0,\nabla\varphi\rangle=0$. Now the equation \eqref{first_eq} becomes
\[
\rho\mathbf{a}_0=\mu|\nabla\varphi|^2\mathbf{a}_0,
\]
which holds true if $\varphi=\varphi_S$ satisfies $c_S^2|\nabla\varphi_S|^2=1$ (cf. \eqref{eikonal0}).

Denote
\[
\begin{split}
(P\mathbf{a})_i=\rho\partial^2_ta_i-&\partial_m(\lambda a_{k;\ell}g^{k\ell}g_{ij}+\mu a_{i;j}+\mu a_{j;i})g^{jm}\\
&+\Gamma_{im}^n(\lambda a_{k;\ell}g^{k\ell}g_{nj}+\mu a_{n;j}+\mu a_{j;n})g^{jm}\\
&+\Gamma_{jm}^n(\lambda a_{k;\ell}g^{k\ell}g_{ni}+\mu a_{n;i}+\mu a_{i;n})g^{jm}.
\end{split}
\]
Notice that $P$ is indeed the elastic wave operator in curvilinear coordinates. Denote also
\[
\begin{split}
(Q\mathbf{a})_i=2\mathrm{i}\rho\partial_ta_i&-\mathrm{i}\left(\partial_i(\lambda a_k\varphi_{;\ell}g^{k\ell})+\lambda a_k\varphi_{;\ell}g^{k\ell}\partial_mg_{ij} g^{jm}+\partial_m(\mu a_i\varphi_{;j}+\mu a_j\varphi_{;i})g^{jm}\right)\\
&-\mathrm{i}\left(\lambda a_{k;\ell}g^{k\ell}\varphi_{;i}+\mu(a_{i;j}+a_{j;i})\varphi_{;m}g^{jm}\right)\\
&+\mathrm{i}\Gamma^n_{im}(\lambda a_k\varphi_{;\ell}g^{k\ell}g_{nj}+\mu a_n\varphi_{;j}+\mu a_j\varphi_{;n})g^{jm}\\
&+\mathrm{i}\Gamma^n_{jm}(\lambda a_k\varphi_{;\ell}g^{k\ell}g_{ni}+\mu a_n\varphi_{;i}+\mu a_i\varphi_{;n})g^{jm},
\end{split}
\]
and
\[
\begin{split}
(R\mathbf{a})_i=-\rho a_i+\lambda a_k\varphi_{;\ell}g^{k\ell}\varphi_{;i}+\mu a_j\varphi_{;m}g^{jm}\varphi_{;i}+\mu\varphi_{;j}\varphi_{;m}g^{jm}a_i).
\end{split}
\]
Setting $\mathbf{a}_{-2}=\mathbf{a}_{-1}=0$, the amplitudes $\mathbf{a}_j$, $j=0,1,2,\cdots$, can be constructed by solving the equations
\[
P\mathbf{a}_{j-2}+Q\mathbf{a}_{j-1}+R\mathbf{a}_{j}=0,\quad j=0,1,2,\cdots.
\]
which is obtained by collecting terms of order $\mathcal{O}(\varrho^{2-j})$ in \eqref{elastic_eq_cova}.
We will use the subprincipal $\mathbf{a}_1$ to recover the density. So we will explicitly construct $\mathbf{a}_0$ and $\mathbf{a}_1$ in the following two sections. The equation $R\mathbf{a}_0=0$ is equivalent to \eqref{first_eq}. For the construction of $\mathbf{a}_0$ and $\mathbf{a}_1$, we also need to use the following equations
\begin{eqnarray}
(R\mathbf{a}_1+Q\mathbf{a}_0)_i=0,\label{second_eq}\\
(R\mathbf{a}_2+Q\mathbf{a}_1+P\mathbf{a}_0)_i=0.\label{third_eq}
\end{eqnarray}

\section{\textit{P}-waves and related results}\label{Pwave}
In this section, we first construct solutions of the form \eqref{solutiondistribution} that represent $P$-waves. We will regain the same results obtained in \cite{rachele2003uniqueness,bhattacharyya2018local}, but the calculation is done in a different manner (under curvilinear coordinates versus Cartesian coordinates). This section is not a part of the proof of the main result of this paper.

For $P$-waves, we take $\varphi=\varphi_P$ such that
\[
|\nabla\varphi|^2=c_P^{-2},
\]
with $\nabla\varphi(x_0)=\xi_0$, $|\xi_0|_{g_P}=1$. We can first calculate the principal term of the amplitude. Denote $J^2=\det(g_{\alpha\beta})$. Assume $\gamma=\gamma_{x_0,\xi_0}$ be the unit-speed geodesic in $(\Omega,g_P)$ issued from $x_0$ in direction $\xi_0$.
\begin{lemma}\label{proofAP}
The principal term of the amplitude can be taken to be
\[
\mathbf{a}_0=A_P\frac{\nabla\varphi}{|\nabla\varphi|}
\]
where $A_P$ is a scalar function of the form
\[
A_P=J^{-1/2}c_P^{1/2}(\lambda+2\mu)^{-1/2}=J^{-1/2}c_P^{-1/2}\rho^{-1/2}.
\]
\end{lemma}
The proofs of the above lemma and all following lemmas are given in the appendix. Observe that one can not get any information about $\rho$ in the interior of $\Omega$ from the value of $\mathbf{a}_0$ at the boundary, which is also pointed out in \cite{rachele2003uniqueness}. To recover the density in the interior, we need to use the subprincipal term of the amplitude $\mathbf{a}_1$, as has been done in \cite{rachele2003uniqueness}.

If we denote the scalar function $B_P=c_P^{3/2}\rho^{1/2}J^{1/2}\langle \mathbf{a}_1,\nabla\varphi\rangle$, we have the following lemma.
\begin{lemma}\label{proofBP}
$B_P=B_P(\gamma(s))$ satisfies the following ordinary differential equation along $\gamma$ in $(\Omega,g_P)$
\[
2\mathrm{i}\frac{\mathrm{d}B_P}{\mathrm{d}s} =M_{ij}\dot{\gamma}^i\dot{\gamma}^j+C,
\]
where $M$ is a second order (covariant) tensor of the form
\begin{align*}
M=&-(c_P^2-4c_S^2)\Delta\log\sqrt{\rho} g_P+\frac{2c_P^2-4c_S^2}{c_P^2}\nabla^2\log\sqrt{\rho}\\
&+\left(c_P^2-4c_S^2+\frac{4c_S^4}{c_P^2-c_S^2}\right)|\nabla\log\sqrt{\rho}|^2g_P+\frac{c_S^2}{c_P^2}\frac{c_P^2-2c_S^2}{c_P^2-c_S^2}\nabla\log\sqrt{\rho}\otimes\nabla\log\sqrt{\rho}\\
&+\nabla\log\sqrt{\rho}\cdot\left(-\nabla c_P^2+\frac{8c_S^2}{c_P^2-c_S^2}\nabla c_S^2\right) g_P+\frac{1}{c_P^2}\nabla\log\sqrt{\rho}\otimes\left(2\nabla c_P^2-\frac{8c_S^2}{c_P^2-c_S^2}\nabla c_S^2\right),
\end{align*}
and $C$ does not depend on $\rho$. 
\end{lemma}

The above lemma has already been proved in \cite{rachele2003uniqueness}. In the appendix we give a different proof by carrying out all the calculations in the ray coordinates.

Recall that $g_P=c_P^{-2}g_E$ is the metric induced by the $P$-wavespeed. If the geodesic $\gamma$ connects two points $a\in\partial\Omega$ to $b\in\partial\Omega$, then
\[
B_P(b)-B_P(a)=\frac{1}{2\mathrm{i}}\int_\gamma (M_{ij}(\gamma(s))\dot{\gamma}^i(s)\dot{\gamma}^j(s)+C)\mathrm{d}s.
\]
By carefully analyzing the reflection of waves (with possible mode conversions) (see, for example, \cite{stefanov2021transmission,caday2021recovery}) at the boundary or using the argument in \cite{rachele2003uniqueness}, one can show that the Dirichlet-to-Neumann map $\Lambda_\rho$ determines $B_P$ on $\partial\Omega$. Similar case is analyzed for acoustic waves in \cite{stefanov2018inverse}. Therefore we can recover the longitudinal geodesic ray transform of $M$. 
\[
I_2M(\gamma)=\int_\gamma M_{ij}(\gamma(s))\dot{\gamma}^i(s)\dot{\gamma}^j(s)\mathrm{d}s.
\]
Here $\gamma$ is an arbitrary geodesic in $(\Omega, g_P)$ connecting two boundary points.

We remark here that the same result was obtained in \cite{rachele2003uniqueness,bhattacharyya2018local}, and then the uniqueness of $\rho$ can be proved assuming, for example, $\lambda=2\mu$ only at isolated points, as well as some geometric conditions. Here one needs to use the injectivity, up to natural obstructions, of the longitudinal geodesic ray transform (cf., for example, \cite{Shara,paternain2015invariant,stefanov2018inverting}) under certain geometric conditions. Indeed, it is proved in \cite{rachele2003uniqueness} that if $I_2M_{\rho_1}=I_2M_{\rho_2}$, then
\[
\begin{split}
&(c_P^2-4c_S^2)\Delta^2(\log\sqrt{\rho_1}-\log\sqrt{\rho_2})\\
&-4\frac{c_P^4-5c_P^2c_S^2+8c_S^4}{c_P^2-c_S^2}\Delta\left[\nabla(\log\sqrt{\rho_1}+\log\sqrt{\rho_2})\cdot \nabla(\log\sqrt{\rho_1}-\log\sqrt{\rho_2})\right]=0.
\end{split}
\]
If $\lambda=2\mu$ only at isolated points, the coefficient $(c_P^2-4c_S^2)$ of the leading order term in the above equation is degenerate only at those isolated points. Then one can use the
boundary determination results  $\partial_\nu^k\rho_1=\partial_\nu^k\rho_2$, $k=0,1,2,3$ on $\partial\Omega$ (cf. \cite{rachele2000boundary}), and the \textit{Unique Continuation Principle} (cf. \cite{protter1960unique}) of elliptic equations to derive $\rho_1=\rho_2$ in $\overline{\Omega}$.

 
\section{\textit{S}-waves and proof of the main theorem}\label{Swave}
In this section, we will use solutions of the form \eqref{solutiondistribution} that represent $S$-waves to derive the main result of this paper. We work under the ray coordinates associated with the metric $g_S=c_S^{-2}\mathrm{d}s^2$. Now we take
\[
|\nabla\varphi|^2=c_S^{-2}
\]
with $\nabla\varphi(x_0)=\xi_0$, $|\xi_0|_{g_S}=1$ and $\langle\mathbf{a}_0,\nabla\varphi\rangle=0$ (or equivalently $a_{03}=0$). Then the identity \eqref{first_eq} is satisfied. Assume $\gamma=\gamma_{x_0,\xi_0}$ be the unit-speed geodesic in $(\Omega,g_S)$ issued from $x_0$ in direction $\xi_0$.

Similar to the $P$-waves, we have

\begin{lemma}\label{proofAS}
The principal term of the amplitude can be taken to be
\[
\mathbf{a}_0(\gamma(s))=A_Sc_S\eta(s)
\]
where $A_S$ is a scalar function of the form
\[
A_S=J^{-1/2}c_S^{-1/2}\rho^{-1/2}
\]
and $\eta$ is a parallel (cotangent) vector field along the geodesic $\gamma$ in $(\Omega,g_S)$ and $\langle \eta,\dot{\gamma}\rangle_{g_S}=0$.
\end{lemma}

Still one can not get any information about $\rho$ in the interior of $\Omega$ from the leading order term of the $S$-wave amplitude at the boundary, as for $P$-waves. We will extract information about $\rho$ in the interior from the subprincipal term $\mathbf{a}_1$.

Denote the scalar function
\[
B_S=J^{1/2}c_S^{-1/2}\rho^{1/2}\langle \mathbf{a}_1,\eta\rangle_{g_S}=J^{1/2}c_S^{3/2}\rho^{1/2}\langle \mathbf{a}_1,\eta\rangle.
\]
Here $\mathbf{a}_1$ and $\eta$ can both be treated as cotangent vectors.
Through tedious calculation we find:
\begin{lemma}\label{proofBS}
$B_S=B_S(\gamma(s))$ satisfies the following ordinary differential equation along $\gamma$
\[
2\mathrm{i}\frac{\mathrm{d}B_S}{\mathrm{d}s} =\langle N,\eta\otimes\eta\rangle_{g_S}+C,
\]
where $N$ is a second (covariant) order tensor of the form
\begin{equation}\label{formofN}
\begin{split}
N=&c_S^2\Delta\log\sqrt{\rho}g_S+2\nabla^2\log\sqrt{\rho}+c_S^2|\nabla\log\sqrt{\rho}|^2g_S+4\frac{c_P^2-2c_S^2}{c_P^2-c_S^2}\nabla\log\sqrt{\rho}\otimes\nabla\log\sqrt{\rho}\\
&+(\nabla\log\sqrt{\rho}\cdot\nabla c_S^2) g_S+\nabla\log\sqrt{\rho}\otimes\left(2\frac{c_P^2-c_S^2}{c_S^4}\nabla c_S^2+c_S^{-2}\nabla c_P^2-\frac{12}{c_P^2-c_S^2}\nabla c_S^2\right),
\end{split}
\end{equation}
and $C$ does not depend on $\rho$. 
\end{lemma}

Here
\[
\langle N,\eta\otimes\eta\rangle_{g_S}=g_S^{ik}g_S^{j\ell}N_{ij}\eta_k\eta_\ell=c_S^4g^{ik}g^{j\ell}N_{ij}\eta_k\eta_\ell.
\]
 If the geodesic $\gamma$ connects two points $a\in\partial\Omega$ to $b\in\partial\Omega$, then
\[
B_S(b)-B_S(a)=\frac{1}{2\mathrm{i}}\int_\gamma (N(\gamma(s)),\eta\otimes\eta\rangle_{g_S}+C)\mathrm{d}s.
\]
 Similar to the discussion for the $P$-waves, we can recover the transverse ray transform of $N$ 
 \[
T_2N(\gamma,\eta)= \int_\gamma \langle N(\gamma(s)),\eta\otimes\eta\rangle_{g_S}\mathrm{d}s
 \]
 from the Dirichlet-to-Neumann map $\Lambda_{\rho}$. Here $\gamma$ is an arbitrary geodesic connecting two boundary points, and $\eta$ could be any parallel vector field along $\gamma$ and orthogonal to $\gamma$ with respect the metric $g_S$ (also with respect to the Euclidean metric $g_E$).  We note here that clearly the transverse ray transform of $N$ is equal to that of the symmetrization of $N$.\\

 Next we switch from ray coordinates to Cartesian coordinates. Under the strictly convex foliation condition, by the invertibility of the transverse ray transform (cf. Proposition \ref{thmglobal20}), we can recover the symmetrization of $N$,
\[
\begin{split}
\mathrm{Sym}(N)=&\Delta\log\sqrt{\rho}I+2\nabla^2\log\sqrt{\rho}+|\nabla\log\sqrt{\rho}|^2I+4\frac{c_P^2-2c_S^2}{c_P^2-c_S^2}\nabla\log\sqrt{\rho}\otimes\nabla\log\sqrt{\rho}\\
&+c_S^{-2}(\nabla\log\sqrt{\rho}\cdot\nabla c_S^2)I+\nabla\log\sqrt{\rho}\otimes^s\left(2\frac{c_P^2-c_S^2}{c_S^4}\nabla c_S^2+c_S^{-2}\nabla c_P^2-\frac{12}{c_P^2-c_S^2}\nabla c_S^2\right).
\end{split}
\]
recalling $c_S^2g_S=g_E$.
Here $I$ is the identity matrix (the Euclidean metric under Cartesian coordinates) and $\otimes^s$ denotes the symmetrized product
\[
u\otimes^s v=\frac{1}{2}(u\otimes v+v\otimes u).
\]
In conclusion, if $\Lambda_{\rho_1}=\Lambda_{\rho_2}$, we have $\mathrm{Sym}(N_1)=\mathrm{Sym}(N_2)$, which can be written as a second order differential equation for $\log\sqrt{\rho_1}-\log\sqrt{\rho_2}$,
\begin{equation}\label{eqrho12}
\begin{split}
&\Delta (\log\sqrt{\rho_1}-\log\sqrt{\rho_2})I+2\nabla^2(\log\sqrt{\rho_1}-\log\sqrt{\rho_2})+(|\nabla\log\sqrt{\rho_1}|^2-|\nabla\log\sqrt{\rho_2}|^2)I\\
&+4\frac{c_P^2-2c_S^2}{c_P^2-c_S^2}\left(\nabla\log\sqrt{\rho_1}\otimes\nabla\log\sqrt{\rho_1}-\nabla\log\sqrt{\rho_2}\otimes\nabla\log\sqrt{\rho_2}\right)\\
&+c_S^{-2}(\nabla(\log\sqrt{\rho_1}-\log\sqrt{\rho_2})\cdot\nabla c_S^2)I\\
&+\left(2\frac{c_P^2-c_S^2}{c_S^4}\nabla c_S^2+c_S^{-2}\nabla c_P^2-\frac{12}{c_P^2-c_S^2}\nabla c_S^2\right)\otimes^s(\nabla\log\sqrt{\rho_1}-\nabla\log\sqrt{\rho_2})\\
=&0\quad \text{in }\Omega.
\end{split}
\end{equation}
Applying divergence twice to the above equation, we obtain a fourth order linear differential equation for $\log\sqrt{\rho_1}-\log\sqrt{\rho_2}$ of the form
\[
3\Delta^2(\log\sqrt{\rho_1}-\log\sqrt{\rho_2})+\sum_{|\alpha|\leq 3}F_\alpha(\rho_1,\rho_2)D^{\alpha}(\log\sqrt{\rho_1}-\log\sqrt{\rho_2})=0,\quad \text{in }\Omega,
\]
by noticing that
\[
\begin{split}
\nabla\cdot\nabla\cdot\left[\Delta (\log\sqrt{\rho_1}-\log\sqrt{\rho_2})I+2\nabla^2(\log\sqrt{\rho_1}-\log\sqrt{\rho_2})\right]=3\Delta^2(\log\sqrt{\rho_1}-\log\sqrt{\rho_2}),
\end{split}
\]
and the remaining terms in \eqref{eqrho12} only involves first order derivatives of $\log\sqrt{\rho_1}-\log\sqrt{\rho_2}$.
 By the \textit{Unique Continuation Principle} (cf. \cite{protter1960unique}) and the boundary determination results  $\partial_\nu^k\rho_1=\partial_\nu^k\rho_2$, $k=0,1,2,3$ on $\partial\Omega$ (cf. \cite{rachele2000boundary}),  we can conclude that
\[
\rho_1=\rho_2\quad\text{in }\overline{\Omega}.
\]
This completes the proof of the main theorem.~\\

\appendix
\section{Proof of lemmas}
In this appendix, we provide all the calculations needed for the proofs of Lemmas \ref{proofAP}-\ref{proofBS}.

\begin{proof}[Proof of Lemma \ref{proofAP}]
In order for \eqref{first_eq} to hold, the principal amplitude $\mathbf{a}_0$ needs to be parallel with $\nabla\varphi$. We take
\[
\mathbf{a}_0=A_P\frac{\nabla\varphi}{|\nabla\varphi|},
\]

Under the ray coordinates associated with $g_P$,
\[
\varphi_{;3}=1,\varphi_{;}^{~~3}=c_P^{-2},\varphi_{;\alpha}=\varphi_{;}^{~~\alpha}=0.
\]
Therefore
\[
a_{03}=c_PA_P\quad a_{0}^{~3}=c_P^{-1}A_P,\quad a_{0\alpha}=a_{0}^{~\alpha}=0,\,\alpha=1,2.
\]

To determine $A_P$, we use the equation \eqref{second_eq} with $i=3$. We calculate
\begin{align*}
&\partial_3(\lambda a_{0k}\varphi_{;\ell}g^{k\ell})+\lambda a_{0k}\varphi_{;\ell}g^{k\ell}\partial_mg_{3j} g^{jm}+\partial_m(\mu a_{03}\varphi_{;j}+\mu a_{0j}\varphi_{;3})g^{jm}\\
&=\partial_3(\lambda c_PA_Pc_P^{-2})+\lambda c_PA_Pc_P^{-2}\partial_3g_{33}g^{33}+\partial_3(\mu c_PA_P+\mu c_PA_P)g^{33}\\
&=\partial_3(\lambda c_P^{-1}A_P)+\lambda c_P^{-3}A_P\partial_3(c_P^2)+\partial_3(2\mu c_PA_P)c_P^{-2}\\
&=\partial_3((\lambda+2\mu) c_PA_P)c_P^{-2},\\
~\\
&\lambda a_{0k;\ell}g^{k\ell}\varphi_{;3}+\mu(a_{03;j}+a_{0j;3})\varphi_{;m}g^{jm}\\
=&\lambda(\partial_3(c_PA_P)c_P^{-2}-c_PA_P\Gamma_{33}^3c_P^{-2}-c_PA_P\Gamma^3_{\alpha\beta}g^{\alpha\beta})+\mu(\partial_3(c_PA_P)-c_PA_P\Gamma_{33}^3+\partial_3(c_PA_P)-c_PA_P\Gamma_{33}^3)c^{-2}_P\\
=&(\lambda+2\mu)\partial_3(c_PA_P)c_P^{-2}-(\lambda+2\mu)A_P(\partial_3c_P)c_P^{-2}+\frac{1}{2}\lambda c_P^{-1}A_P\frac{\partial g_{\alpha\beta}}{\partial\tau}g^{\alpha\beta},\\
~\\
&\Gamma^n_{3m}(\lambda a_{0k}\varphi_{;\ell}g^{k\ell}g_{nj}+\mu a_{0n}\varphi_{;j}+\mu a_{0j}\varphi_{;n})g^{jm}\\
=&\Gamma_{3n}^n\lambda c_PA_Pg^{33}+\Gamma_{33}^3(\mu c_PA_P+\mu c_PA_P)g^{33}\\
=&\frac{1}{2}g^{\alpha\beta}\frac{\partial g_{\alpha\beta}}{\partial\tau}\lambda c_P^{-1}A_P+c_P^{-2}\partial_3c_P(\lambda+2\mu)A_P,\\
~\\
&\Gamma^n_{jm}(\lambda a_{0k}\varphi_{;\ell}g^{k\ell}g_{n3}+\mu a_{0n}\varphi_{;3}+\mu a_{03}\varphi_{;n})g^{jm}\\
=&\Gamma_{jm}^3(\lambda a_{03}g^{33}g_{33}+\mu a_{03}+\mu a_{03})g^{jm}\\
=&(c_P^{-3}\partial_3c_P-\frac{1}{2}c_P^{-2}\frac{\partial g_{\alpha\beta}}{\partial\tau}g^{\alpha\beta})(\lambda+2\mu)c_PA_P.
\end{align*}


Summing up the above calculations, one obtains
\[
\begin{split}
&-(\lambda a_{13}c_P^{-2}+2\mu a_{13}c_P^{-2}-\rho a_{13})\\
&+\mathrm{i}c_P^{-1}\left(2(\lambda+2\mu)(\partial_3-\partial_t)A_P-(\lambda+2\mu)c_P^{-1}(\partial_3c_P)A_P+\partial_3(\lambda+2\mu)A_P+\frac{1}{2}(\lambda+2\mu)\frac{\partial g_{\alpha\beta}}{\partial\tau}g^{\alpha\beta}\right)=0.
\end{split}
\]

We can take $A_P$ to be independent of $t$.
Using the equality
\[
g^{\alpha\beta}\frac{\partial g_{\alpha\beta}}{\partial\tau}=\frac{2}{J}\frac{\partial J}{\partial\tau},
\]
we end up with
\[
J(\lambda+2\mu)\partial_\tau A_P-\frac{1}{2}J(\lambda+2\mu)c_P^{-1}(\partial_\tau c_P)A_P+\frac{1}{2}J\partial_\tau(\lambda+2\mu)A_P+\frac{1}{2}(\partial_\tau J)(\lambda+2\mu)A_P=0.
\]
The above equation can be considered as a linear ODE along a single geodesic $\gamma$.
Note that it can be written as
\[
\frac{\partial}{\partial \tau}(A_PJ^{1/2}c_P^{-1/2}(\lambda+2\mu)^{1/2})=0.
\]
Solving the above equation, one gets
\[
A_P=C_0J^{-1/2}c_P^{1/2}(\lambda+2\mu)^{-1/2}=C_0J^{-1/2}c_P^{-1/2}\rho^{-1/2}
\]
with an arbitrary constant $C_0$. This completes the proof of the lemma. Since we need to use nontrivial solutions, one can just take $C_0=1$.
\end{proof}
~\\

\begin{proof}[Proof of Lemma \ref{proofBP}]
In order the construct the subprincipal term of the amplitude, we first use equation \eqref{second_eq} with $i=\alpha$. For this, we calculate
\begin{align*}
&\partial_\alpha(\lambda a_{0k}\varphi_{;\ell}g^{k\ell})+\lambda a_{0k}\varphi_{;\ell}g^{k\ell}\partial_mg_{\alpha j} g^{jm}+\partial_m(\mu a_{0\alpha}\varphi_{;j}+\mu a_{0j}\varphi_{;\alpha})g^{jm}\\
=&\partial_\alpha(\lambda c_PA_P c_P^{-2})+\lambda c_PA_Pc_P^{-2}(\partial_\varepsilon g_{\alpha\beta})g^{\beta\varepsilon},\\
~\\
&\lambda a_{0k;\ell}g^{k\ell}\varphi_{;\alpha}+\mu(a_{0\alpha;j}+a_{0j;\alpha})\varphi_{;m}g^{jm}\\
=&\mu(-\Gamma_{\alpha 3}^3a_{03}+\partial_\alpha (a_{03})-\Gamma_{\alpha 3}^3a_{03})g^{33}\\
=&\mu(\partial_\alpha(c_PA_P) -2c_P^{-1}(\partial_\alpha c_P)c_PA_P)c_P^{-2},\\
~\\
&\Gamma^n_{\alpha m}(\lambda a_{0k}\varphi_{;\ell}g^{k\ell}g_{nj}+\mu a_{0n}\varphi_{;j}+\mu a_{0j}\varphi_{;n})g^{jm}\\
=&\Gamma_{\alpha n}^n(\lambda c_PA_Pc_P^{-2})+\Gamma_{\alpha 3}^3 g^{33}(2\mu c_PA_P)\\
=&\Gamma_{\alpha\beta}^\beta\lambda c_P^{-1}A_P+\Gamma_{\alpha 3}^3 (\lambda+2\mu) c_P^{-1}A_P,\\
~\\
&\Gamma^n_{jm}(\lambda a_{0k}\varphi_{;\ell}g^{k\ell}g_{n\alpha}+\mu a_{0n}\varphi_{;\alpha}+\mu a_{0\alpha}\varphi_{;n})g^{jm}\\
=&\Gamma_{jm}^\beta(\lambda a_{03}g^{33}g_{\alpha\beta})g^{jm}\\
=&\Gamma^\beta_{33}\lambda c_P^{-1}A_Pg_{\alpha\beta}c_P^{-2}+\Gamma^\beta_{\varepsilon\delta}\lambda c_P^{-1}A_P g_{\alpha\beta}g^{\varepsilon\delta}.
\end{align*}

Using the identity
\begin{equation}\label{Christofel_identity}
(\partial_\varepsilon g_{\alpha\beta})g^{\beta\gamma}-\Gamma_{\alpha\beta}^\beta-\Gamma^\beta_{\varepsilon\delta}g_{\alpha\beta}g^{\varepsilon\delta}=0,
\end{equation}
we obtain
\[
\begin{split}
&-\mu c_P^{-2}a_{1\alpha}+\rho a_{1\alpha}\\
&+\mathrm{i}(\partial_\alpha(\lambda A_P )+\mu(\partial_\alpha(c_PA_P) -2(\partial_\alpha c_P)A_P)c_P^{-2}-\Gamma_{\alpha 3}^3 (\lambda+2\mu) c_P^{-1}A_P-\Gamma^\beta_{33}\lambda c_P^{-3}A_Pg_{\alpha\beta})=0,
\end{split}
\]
which simplifies to
\[
-\mu c_P^{-2}a_{1\alpha}+\rho a_{1\alpha}+\mathrm{i}((\lambda+\mu)c_P^{-1}\partial_\alpha A_P+(\partial_\alpha\lambda)c_P^{-1}A_P-c_P^{-2}\partial_\alpha c_P(3\mu+\lambda)A_P)=0.
\]

So, for $\alpha=1,2$,
\begin{equation}\label{eq_a1alpha}
\begin{split}
a_{1\alpha}=&\frac{\mathrm{i}}{c_S^{2}c_P^{-2}- 1}\left(\frac{\lambda+\mu}{\rho}c_P^{-1}\partial_\alpha A_P+\frac{(\partial_\alpha\lambda)}{\rho}c_P^{-1}A_P-c_P^{-2}\partial_\alpha c_P\frac{3\mu+\lambda}{\rho}A_P\right)\\
=&\frac{\mathrm{i}}{c_S^{2}c_P^{-2}- 1}\Big((c_P^2-c_S^2)c_P^{-1}J^{-1/2}c_P^{-1/2}(-\frac{1}{2}\rho^{-3/2}\partial_\alpha\rho)+(c_P^2-c_S^2)c_P^{-1}\partial_\alpha(J^{-1/2}c_P^{-1/2})\rho^{-1/2}\\
&+\rho^{-3/2}(c_P^2-2c_S^2)(\partial_\alpha\rho)c_P^{-1}J^{-1/2}c_P^{-1/2}+\partial_\alpha(c_P^2-2c_S^2)c_P^{-1}J^{-1/2}c_{P}^{-1/2}\rho^{-1/2}\\
&-c_P^{-2}(\partial_\alpha c_P)(c_P^2+c_S^2)J^{-1/2}c_P^{-1/2}\rho^{-1/2}\Big)\\
=&\frac{\mathrm{i}}{c_S^{2}c_P^{-2}- 1}\Big(\frac{1}{2}\rho^{-3/2}(c_P^2-3c_S^2)(\partial_\alpha\rho)c_P^{-3/2}J^{-1/2}+(c_P^2-c_S^2)c_P^{-3/2}\partial_\alpha(J^{-1/2})\rho^{-1/2}\\
&\quad\quad\quad\quad+c_P^{-5/2}\left(\frac{1}{2}c_P^2\partial_\alpha c_P-\frac{1}{2}c_S^2\partial_\alpha c_P-4c_Sc_P\partial_\alpha c_S\right)J^{-1/2}\rho^{-1/2}\Big)\\
=&\frac{\mathrm{i}}{(c_S^{2}c_P^{-2}- 1)\sqrt{\rho}}\Big((c_P^2-3c_S^2)(\partial_\alpha\log\sqrt{\rho})c_P^{-3/2}J^{-1/2}+(c_P^2-c_S^2)c_P^{-3/2}\partial_\alpha(J^{-1/2})\\
&\quad\quad\quad\quad+c_P^{-5/2}\left(\frac{1}{2}c_P^2\partial_\alpha c_P-\frac{1}{2}c_S^2\partial_\alpha c_P-4c_Sc_P\partial_\alpha c_S\right)J^{-1/2}\Big)\\
=&\frac{\mathrm{i}}{(c_S^{2}- c_P^2)\sqrt{\rho}}\Big((c_P^2-3c_S^2)(\partial_\alpha\log\sqrt{\rho})c_P^{1/2}J^{-1/2}+(c_P^2-c_S^2)c_P^{1/2}\partial_\alpha(J^{-1/2})\\
&\quad\quad\quad\quad+c_P^{-1/2}\left(\frac{1}{2}c_P^2\partial_\alpha c_P-\frac{1}{2}c_S^2\partial_\alpha c_P-4c_Sc_P\partial_\alpha c_S\right)J^{-1/2}\Big).
\end{split}
\end{equation}

Now we have $a_{11}$ and $a_{12}$. In order the obtain $a_{13}$,  we subsequently use equation \eqref{third_eq} with $i=3$, that is, 

\begin{align*}
&2\mathrm{i}\rho\partial_ta_{13}\\
=&\mathrm{i}\left(\partial_3(\lambda a_{1k}\varphi_{;\ell}g^{k\ell})+\lambda a_{1k}\varphi_{;\ell}g^{k\ell}\partial_mg_{3j} g^{jm}+\partial_m(\mu a_{13}\varphi_{;j}+\mu a_{1j}\varphi_{;3})g^{jm}\right)\\
&\mathrm{i}\left(\lambda a_{1k;\ell}g^{k\ell}\varphi_{;3}+\mu(a_{13;j}+a_{1j;3})\varphi_{;m}g^{jm}\right)\\
&-\mathrm{i}\Gamma^n_{3m}(\lambda a_{1k}\varphi_{;\ell}g^{k\ell}g_{nj}+\mu a_{1n}\varphi_{;j}+\mu a_{1j}\varphi_{;n})g^{jm}\\
&-\mathrm{i}\Gamma^n_{jm}(\lambda a_{1k}\varphi_{;\ell}g^{k\ell}g_{n3}+\mu a_{1n}\varphi_{;3}+\mu a_{13}\varphi_{;n})g^{jm}\\
&+\partial_m(\lambda a_{0k;\ell}g^{k\ell}g_{3j}+\mu a_{03;j}+\mu a_{0j;3})g^{jm}\\
&-\Gamma_{3m}^n(\lambda a_{0k;\ell}g^{k\ell}g_{nj}+\mu a_{0n;j}+\mu a_{0j;n})g^{jm}\\
&-\Gamma_{jm}^n(\lambda a_{0k;\ell}g^{k\ell}g_{n3}+\mu a_{0n;3}+\mu a_{3;n})g^{jm}.
\end{align*}
Then we calculate
\begin{align*}
\partial_\beta a_{1\alpha}=&\frac{\mathrm{i}}{(c_S^{2}- c_P^{2})\sqrt{\rho}}\Bigg((c_P^2-3c_S^2)(\partial_\beta\partial_\alpha\log\sqrt{\rho})c_P^{1/2}J^{-1/2}+\partial_\beta((c_P^2-3c_S^2)c_P^{1/2}J^{-1/2})\partial_\alpha\log\sqrt{\rho}\\
&\quad\quad\quad\quad\quad-\frac{\partial_\beta c_S^2-\partial_\beta c_P^2}{(c_S^2-c_P^2)}(c_P^2-3c_S^2)(\partial_\alpha\log\sqrt{\rho})c_P^{1/2}J^{-1/2}\\
&\quad\quad\quad\quad\quad+(c_P^2-3c_S^2)(\partial_\alpha\log\sqrt{\rho})c_P^{1/2}J^{-1/2}(-\partial_\beta\log\sqrt{\rho})\\
&\quad\quad\quad\quad\quad+(c_P^2-c_S^2)c_P^{1/2}\partial_\alpha(J^{-1/2})(-\partial_\beta\log\sqrt{\rho})\\
&\quad\quad\quad\quad\quad+c_P^{-1/2}\left(\frac{1}{2}c_P^2\partial_\alpha c_P-\frac{1}{2}c_S^2\partial_\alpha c_P-4c_Sc_P\partial_\alpha c_S\right)J^{-1/2}(-\partial_\beta\log\sqrt{\rho})+C\Bigg).
\end{align*}

With the above expression, we get
\begin{align*}
&\partial_3(\lambda a_{1k}\varphi_{;\ell}g^{k\ell})+\lambda a_{1k}\varphi_{;\ell}g^{k\ell}\partial_mg_{3j} g^{jm}+\partial_m(\mu a_{13}\varphi_{;j}+\mu a_{1j}\varphi_{;3})g^{jm}\\
=&\partial_3(\lambda a_{13}g^{33})+\lambda a_{13}g^{33}\partial_3(g_{33})g^{33}+\partial_3(\mu a_{13})g^{33}+\partial_m(\mu a_{1j})g^{jm}\\
=&\partial_3(\lambda a_{13}c_P^{-2})+\lambda a_{13} c_P^{-4}\partial_3(c_P^2)+2\partial_3(\mu a_{13})c_P^{-2}+\partial_\alpha(\mu a_{1\beta})g^{\alpha\beta}\\
=&\partial_3((\lambda+2\mu) a_{13})c_P^{-2}+\partial_\alpha(\mu a_{1\beta})g^{\alpha\beta},\\
~\\
&\lambda a_{1k;\ell}g^{k\ell}\varphi_{;3}+\mu(a_{13;j}+a_{1j;3})\varphi_{;m}g^{jm}\\
=&\lambda(a_{13;3}g^{33}+a_{1\alpha;\beta}g^{\alpha\beta})+2\mu a_{13;3}g^{33}\\
=&\lambda(\partial_3a_{13}-\Gamma^j_{33}a_{1j})c_P^{-2}+\lambda(\partial_\beta a_{1\alpha}-\Gamma^j_{\alpha\beta}a_{1j})g^{\alpha\beta}+2\mu (\partial_3a_{13}-\Gamma_{33}^ja_{1j})c_P^{-2}\\
=&(\lambda+2\mu)c_P^{-2}\partial_3a_{13}-(\lambda+2\mu)c_P^{-1}(\partial_3c_P)c_P^{-2}a_{13}-(\lambda+2\mu)(\Gamma_{33}^\alpha a_{1\alpha})c_P^{-2}+\lambda g^{\alpha\beta}\partial_\beta a_{1\alpha}\\
&+\frac{1}{2}\lambda c_P^{-2}a_{13}\frac{\partial g_{\alpha\beta}}{\partial\tau}g^{\alpha\beta}-\lambda g^{\alpha\beta}\Gamma_{\alpha\beta}^\varepsilon a_{1\varepsilon},\\
~\\
&\Gamma^n_{3m}(\lambda a_{1k}\varphi_{;\ell}g^{k\ell}g_{nj}+\mu a_{1n}\varphi_{;j}+\mu a_{1j}\varphi_{;n})g^{jm}\\
=&\Gamma_{3n}^n\lambda a_{13}g^{33}+\Gamma_{33}^n\mu a_{1n}g^{33}+\Gamma_{3m}^3\mu a_{1j}g^{jm}\\
=&\frac{1}{2}g^{\alpha\beta}\frac{\partial g_{\alpha\beta}}{\partial\tau}\lambda c_P^{-2}a_{13}+c_P^{-3}\partial_3c_P(\lambda+2\mu)a_{13}+\Gamma_{33}^\alpha\mu a_{1\alpha}c_P^{-2}+\Gamma_{3\beta}^3\mu a_{1\alpha}g^{\alpha\beta},\\
~\\
&\Gamma^n_{jm}(\lambda a_{1k}\varphi_{;\ell}g^{k\ell}g_{n3}+\mu a_{1n}\varphi_{;3}+\mu a_{13}\varphi_{;n})g^{jm}\\
=&\Gamma_{jm}^3(\lambda a_{13}g^{33}g_{33}+\mu a_{13}+\mu a_{13})g^{jm}+\Gamma_{jm}^\alpha(\mu a_{1\alpha})g^{jm}\\
=&(c_P^{-3}\partial_3c_P-\frac{1}{2}c_P^{-2}\frac{\partial g_{\alpha\beta}}{\partial\tau}g^{\alpha\beta})(\lambda+2\mu)a_{13}+\Gamma^\alpha_{\beta\varepsilon}\mu a_{1\alpha}g^{\beta\varepsilon}+\Gamma^\alpha_{33}\mu a_{1\alpha}c_P^{-2}.
\end{align*}

We also calculate
\begin{align*}
&\partial_m(\lambda a_{0k;\ell}g^{k\ell}g_{3j}+\mu a_{03;j}+\mu a_{0j;3})g^{jm}\\
=&\partial_3(\lambda(\partial_3a_{03}-\Gamma^j_{33}a_j)g^{33}g_{33}+\lambda(-\Gamma^j_{\alpha\beta}a_j)g^{\alpha\beta}g_{33})g^{33}\\
&\quad\quad\quad+\partial_m(\mu\partial_j a_{03}-\mu\Gamma^i_{3j}a_{0i})g^{jm}+\partial_m(\mu \partial_3 a_{0j}-\mu\Gamma_{3j}^ia_{0i})g^{jm}\\
=&\partial_3(\lambda\partial_3(c_PA_P)-\lambda\Gamma^3_{33}c_PA_P-\lambda\Gamma^3_{\alpha\beta}c_PA_Pg^{\alpha\beta}c_P^2)c_P^{-2}\\
&+g^{jm}\partial_m(\mu\partial_j(c_PA_P))+g^{33}\partial_3(\mu\partial_3(c_PA_P))-2g^{jm}\partial_m(\mu\Gamma^3_{3j}c_PA_P)\\
=&\partial_3\left(\lambda\partial_3(c_PA_P)-\lambda(\partial_3c_P)A_P+\frac{1}{2}\lambda c_Pg^{\alpha\beta}\frac{\partial g_{\alpha\beta}}{\partial \tau}A_P\right)c_P^{-2}+2c_P^{-2}\partial_3(\mu\partial_3(c_PA_P))+g^{\alpha\beta}\partial_\beta(\mu\partial_\alpha(c_PA_P))\\
&-2c_P^{-2}\partial_3(\mu (\partial_3c_P)A_P)-2g^{\alpha\beta}\partial_\alpha(\mu(\partial_\beta c_P)A_P)\\
=&\partial_3((\lambda+2\mu)c_P\partial_3A_P)c_P^{-2}+\frac{1}{2}c_P^{-1}\lambda g^{\alpha\beta}\frac{\partial g_{\alpha\beta}}{\partial\tau}\partial_3A_P+g^{\alpha\beta}\partial_\beta(\mu\partial_\alpha(c_PA_P))+\frac{1}{2}\partial_3\lambda c_P^{-1}g^{\alpha\beta}\frac{\partial g_{\alpha\beta}}{\partial\tau}A_P\\
&-2g^{\alpha\beta}\mu(\partial_\beta c_P)\partial_\alpha A_P-2g^{\alpha\beta}\partial_\alpha\mu\partial_\beta c_PA_P+\sqrt{\rho}\times\text{terms that do not depend on $\rho$},\\
~\\
&\Gamma_{3m}^n(\lambda a_{0k;\ell}g^{k\ell}g_{nj}+\mu a_{0n;j}+\mu a_{0j;n})g^{jm}\\
=&\Gamma_{3n}^n(\lambda(\partial_3 a_{03}-\Gamma^j_{33}a_j)g^{33}-\lambda\Gamma^j_{\alpha\beta}a_jg^{\alpha\beta})+\Gamma_{3m}^n(\mu \partial_j a_{0n}+\mu\partial_n a_{0j}-2\mu\Gamma^i_{nj}a_{0i})g^{jm}\\
=&\Gamma_{3n}^n(\lambda\partial_3 (c_PA_P)g^{33})+\Gamma^3_{3m}(\mu\partial_j(c_PA_P))g^{jm}+\Gamma^n_{33}(\mu\partial_n(c_PA_P))c_P^{-2}\\
&+\sqrt{\rho}\times\text{terms that do not depend on $\rho$}\\
=&\Gamma_{33}^3(\lambda+2\mu)c_P\partial_3A_P c_P^{-2}+\Gamma_{3\alpha}^\alpha \lambda c_P\partial_3A_P c_P^{-2}+\Gamma_{3\alpha}^3g^{\alpha\beta}\mu c_P\partial_\beta A_P+\Gamma_{33}^\alpha \mu c_P\partial_\alpha A_P\\
&+\sqrt{\rho}\times\text{terms that do not depend on $\rho$}\\
=&c_P^{-2}(\partial_3c_P)(\lambda+2\mu)(\partial_3A_P)+\frac{1}{2}g^{\delta\varepsilon}\frac{\partial g_{\delta\varepsilon}}{\partial\tau}\lambda c_P^{-1}(\partial_3A_P)+\sqrt{\rho}\times\text{terms that do not depend on $\rho$},\\
~\\
&\Gamma_{jm}^n(\lambda a_{0k;\ell}g^{k\ell}g_{n3}+\mu a_{0n;3}+\mu a_{03;n})g^{jm}\\
=&g^{jm}\Gamma^3_{jm}(\lambda+2\mu)\partial_3 (c_PA_P)+g^{jm}\Gamma^\alpha_{jm}\mu\partial_\alpha(c_PA_P)+\sqrt{\rho}\times\text{terms that do not depend on $\rho$}\\
=&c_P^{-2}c_P^{-1}(\partial_3c_P)(\lambda+2\mu)c_P\partial_3A_P+g^{\alpha\beta}(-\frac{1}{2}c_P^{-2}\frac{\partial g_{\alpha\beta}}{\partial\tau})(\lambda+2\mu) c_P\partial_3A_P\\
&+c_P^{-2}(-c_Pg^{\alpha\beta}(\partial_\beta c_P))\mu c_P\partial_\alpha A_P+g^{\varepsilon\delta}\Gamma^\alpha_{\varepsilon\delta}\mu c_P\partial_\alpha A_P+\sqrt{\rho}\times\text{terms that do not depend on $\rho$}\\
=&c_P^{-2}(\partial_3c_P)(\lambda+2\mu)(\partial_3A_P)-\frac{1}{2}g^{\alpha\beta}\frac{\partial g_{\alpha\beta}}{\partial\tau}(\lambda+2\mu) c_P^{-1}(\partial_3A_P)\\
&-g^{\alpha\beta}(\partial_\beta c_P)\mu (\partial_\alpha A_P)+g^{\varepsilon\delta}\Gamma^\alpha_{\varepsilon\delta}\mu c_P(\partial_\alpha A_P)+\sqrt{\rho}\times\text{terms that do not depend on $\rho$}.
\end{align*}

We end up with an equation for $a_{13}$, which is taken to be independent of $t$,
\begin{equation}\label{transport_eq_amplitidue1P}
\begin{split}
&2\mathrm{i}(\lambda+2\mu)c_P^{-2}\partial_3a_{13}-3\mathrm{i}(\lambda+2\mu)c_P^{-3}\partial_3 c_P a_{13}+\mathrm{i}\partial_3(\lambda+2\mu)c_P^{-2} a_{13}+\frac{1}{2}\mathrm{i}c_P^{-2}g^{\alpha\beta}\frac{g_{\alpha\beta}}{\partial\tau}(\lambda+2\mu)a_{13}\\
&+\mathrm{i}H_1+H_2=0,
\end{split}
\end{equation}
where
\[
H_1=(\lambda+\mu) g^{\alpha\beta}\partial_\alpha a_{1\beta}+g^{\alpha\beta}\partial_\alpha\mu a_{1\beta}-(\lambda+4\mu)c_P^{-2}\Gamma^\alpha_{33}a_{1\alpha}-(\lambda+\mu)g^{\alpha\beta}\Gamma^\varepsilon_{\alpha\beta}a_{1\varepsilon}-\mu c_P^{-2}\Gamma^\alpha_{33}a_{1\alpha},
\]
and
\begin{align*}
H_2=&(\lambda+2\mu)c_P^{-1}\partial^2_{33}A_P+\partial_3(\lambda+2\mu)c_P^{-1}\partial_3A_P-c_P^{-2}(\partial_3c_P)(\lambda+2\mu)(\partial_3A_P)+g^{\alpha\beta}\partial_\beta(\mu\partial_\alpha(c_PA_P))\\
&+\frac{1}{2}\partial_3\lambda c_P^{-1}g^{\alpha\beta}\frac{\partial g_{\alpha\beta}}{\partial \tau}A_P-2g^{\alpha\beta}\partial_\alpha\mu\partial_\beta c_PA_P\\
&-g^{\alpha\beta}\mu \partial_\beta c_P\partial_\alpha A_P
+\frac{1}{2}g^{\alpha\beta}\frac{\partial g_{\alpha\beta}}{\partial\tau}(\lambda+2\mu) c_P^{-1}(\partial_3A_P)-g^{\varepsilon\delta}\Gamma^\alpha_{\varepsilon\delta}\mu c_P(\partial_\alpha A_P)\\
&+\sqrt{\rho}\times\text{terms that do not depend on $\rho$}.
\end{align*}
Inserting the expression of $a_{1\alpha}$ \eqref{eq_a1alpha}, we obtain
\begin{align*}
H_1=\frac{\mathrm{i}\sqrt{\rho}}{(c_S^2-c_P^2)}\Bigg( (c_P^2-3c_S^2)g^{\alpha\beta}\partial_\alpha\partial_\beta\log\sqrt{\rho}c_P^{1/2}J^{-1/2}+g^{\alpha\beta}\partial_\beta\left((c_P^2-3c_S^2)c_P^{1/2}J^{-1/2})\right)\partial_\alpha\log\sqrt{\rho}\\
-\frac{\partial_\beta c_S^2-\partial_\beta c_P^2}{(c_S^2-c_P^2)}(c_P^2-3c_S^2)(\partial_\alpha\log\sqrt{\rho})c_P^{1/2}J^{-1/2}\\
-(c_P^2-3c_S^2)g^{\alpha\beta}(\partial_\alpha\log\sqrt{\rho})(\partial_\beta\log\sqrt{\rho}) c_P^{1/2}J^{-1/2}-(c_P^2-c_S^2)c_P^{1/2}g^{\alpha\beta}\partial_\alpha(J^{-1/2})\partial_\beta\log\sqrt{\rho}\\
-\frac{1}{2}c_P^{-1/2}(c_P^2-c_S^2)g^{\alpha\beta}(\partial_\alpha c_P)(\partial_\beta\log\sqrt{\rho})J^{-1/2}+4c_S c_P^{1/2}g^{\alpha\beta}(\partial_\alpha c_S)J^{-1/2}(\partial_\beta\log\sqrt{\rho})\Bigg)\\
+g^{\alpha\beta}(2c_S^2\rho\partial_\alpha\log\sqrt{\rho}+2c_S\partial_\alpha c_S\rho)\frac{\mathrm{i}}{(c_S^2-c_P^2)\sqrt{\rho}}\left((c_P^2-3c_S^2)(\partial_\beta\log\sqrt{\rho})c_P^{1/2}J^{-1/2}\right)\\
+2g^{\alpha\beta}c_S^2\rho\partial_\alpha\log\sqrt{\rho}\frac{\mathrm{i}}{(c_S^{2} -c_P^2)\sqrt{\rho}}(c_P^2-c_S^2)c_P^{1/2}\partial_\beta(J^{-1/2})\\
+2g^{\alpha\beta}c_S^2\rho\partial_\beta\log\sqrt{\rho}\frac{\mathrm{i}}{(c_S^{2}- c_P^2)\sqrt{\rho}}c_P^{-1/2}\left(\frac{1}{2}c_P^2\partial_\alpha c_P-\frac{1}{2}c_S^2\partial_\alpha c_P-4c_Sc_P\partial_\alpha c_S\right)J^{-1/2}\\
-(\lambda+4\mu)(-c_Pg^{\alpha\beta}\partial_\beta c_P)c_P^{-2}\frac{\mathrm{i}}{(c_S^2-c_P^2)\sqrt{\rho}}\left((c_P^2-3c_S^2)(\partial_\alpha\log\sqrt{\rho})c_P^{1/2}J^{-1/2}\right)\\
-(\lambda+\mu)g^{\alpha\beta}\Gamma^\varepsilon_{\alpha\beta}\frac{\mathrm{i}}{(c_S^2-c_P^2)\sqrt{\rho}}\left((c_P^2-3c_S^2)(\partial_\varepsilon\log\sqrt{\rho})c_P^{1/2}J^{-1/2}\right)\\
-c_P^{-1}\partial_\alpha c_P\mu g^{\alpha\beta}\frac{\mathrm{i}}{(c_S^2-c_P^2)\sqrt{\rho}}\left((c_P^2-3c_S^2)(\partial_\beta\log\sqrt{\rho})c_P^{1/2}J^{-1/2}\right)\\
+\sqrt{\rho}\times\text{terms that do not depend on $\rho$}.
\end{align*}
The above formula can be simplified to
\[
\begin{split}
H_1=&-\mathrm{i}\sqrt{\rho}(c_P^2-3c_S^2)c_P^{1/2}J^{-1/2}g^{\alpha\beta}\partial_\alpha\partial_\beta\log\sqrt{\rho}\\
&+\mathrm{i}\left(c_P^2-3c_S^2-\frac{2(c_P^2-3c_S^2)}{c_P^2-c_S^2}\right)\sqrt{\rho}c_S^2 c_P^{1/2}J^{-1/2}g^{\alpha\beta}\partial_\alpha\log\sqrt{\rho}\partial_\beta\log\sqrt{\rho}\\
&+\mathrm{i}g^{\alpha\beta}\Gamma^\varepsilon_{\alpha\beta}(c_P^2-3c_S^2)c_P^{1/2}J^{-1/2}\partial_\varepsilon\log\sqrt{\rho}\\
&-\mathrm{i}\sqrt{\rho}(c_P^2-3c_S^2)c_P^{-1/2}J^{-1/2}g^{\alpha\beta}\partial_\alpha\log\sqrt{\rho}\partial_\beta c_P\\
&+\mathrm{i}\sqrt{\rho}\left(2c_S-\frac{4c_S(c_P^2-5c_S^2)}{c_P^2-c_S^2}\right)c_P^{-1/2}J^{-1/2}g^{\alpha\beta}\partial_\alpha\log\sqrt{\rho}\partial_\beta c_S\\
&+\sqrt{\rho}\times\text{terms that do not depend on $\rho$}.
\end{split}
\]

Further, using the identities
\[
\partial^2_{33}\rho^{-1/2}=-\rho^{-1/2}\partial^2_{33}\log\sqrt{\rho}+\rho^{-1/2}(\partial_3\log\sqrt{\rho})^2,
\]
\[
\partial_\alpha\partial_\beta\rho^{-1/2}=-\rho^{-1/2}\partial_\alpha\partial_\beta\log\sqrt{\rho}+\rho^{-1/2}(\partial_\alpha\log\sqrt{\rho})(\partial_\beta\log\sqrt{\rho}),
\]
we have
\begin{align*}
H_2=&(\lambda+2\mu)c_P^{-1}\left(2(\partial_3J^{-1/2})c_P^{-1/2}\partial_3\rho^{-1/2}+2J^{-1/2}(\partial_3c_P^{-1/2})\partial_3\rho^{-1/2}+J^{-1/2}c_P^{-1/2}\partial^2_{33}\rho^{-1/2}\right)\\
&+\left(2(\lambda+2\mu)\partial_3\log\sqrt{\rho}+\partial_3(c_P^2)\rho\right)c_P^{-1}J^{-1/2}c_P^{-1/2}(-\frac{1}{2}\rho^{-1/2}\partial_3\log\sqrt{\rho})\\
&+2(\lambda+2\mu)\partial_3\log\sqrt{\rho}c_P^{-1}\partial_3(c_P^{-1/2}J^{-1/2})\rho^{-1/2}\\
&-c_P^{-2}(\partial_3c_P)(\lambda+2\mu)J^{-1/2}c_P^{-1/2}(-\rho^{-1/2}\partial_3\log\sqrt{\rho})\\
&+g^{\alpha\beta}(\partial_\beta\mu)\left(\partial_\alpha J^{-1/2}c_P^{1/2}\rho^{-1/2}+J^{-1/2}\partial c_P^{1/2}\rho^{-1/2}+J^{-1/2}c_P^{1/2}\partial_\alpha\rho^{-1/2}\right)\\
&+\mu g^{\alpha\beta}\left(2\partial_\alpha J^{-1/2}c_P^{1/2}\partial_\beta\rho^{-1/2}+2J^{-1/2}\partial_\alpha c_P^{1/2}\partial_\beta\rho^{-1/2}+J^{-1/2}c_P^{1/2}\partial_\alpha\partial_\beta\rho^{-1/2}\right)\\
&+(c_P^2-2c_S^2)\rho\partial_3\log\sqrt{\rho} c_P^{-1}g^{\alpha\beta}\frac{\partial g^{\alpha\beta}}{\partial\tau}c_P^{-1/2}J^{-1/2}\rho^{-1/2}\\
&-4g^{\alpha\beta}\mu \partial_\beta c_PJ^{-1/2}c_P^{-1/2}\partial_\alpha\log\sqrt{\rho}\\
&-g^{\alpha\beta}(\partial_\beta c_P)\mu J^{-1/2}c_P^{-1/2}(-\rho^{-1/2}\partial_\alpha\log\sqrt{\rho})\\
&+\frac{1}{2}g^{\alpha\beta}\frac{\partial g_{\alpha\beta}}{\partial\tau}(\lambda+2\mu) c_P^{-1}J^{-1/2}c_P^{-1/2}\left(-\rho^{-1/2}\partial_3\log\sqrt{\rho}\right)-g^{\varepsilon\delta}\Gamma^\alpha_{\varepsilon\delta}\mu c_PJ^{-1/2}c_P^{-1/2}(-\rho^{-1/2}\partial_\alpha\log\sqrt{\rho})\\
&+\sqrt{\rho}\times\text{terms that do not depend on $\rho$}\\
=&2(\lambda+2\mu)c_P^{-1}(\partial_3J^{-1/2})c_P^{-1/2}(-\rho^{-1/2}\partial_3\log\sqrt{\rho})+(\lambda+2\mu)c_P^{-1}J^{-1/2}c_P^{-3/2}(\partial_3c_P)(-\rho^{-1/2}\partial_3\log\sqrt{\rho})\\
&+(\lambda+2\mu)c_P^{-1}J^{-1/2}c_P^{-1/2}\rho^{-1/2}(-\partial^2_{33}\log\sqrt{\rho}+(\partial_3\log\sqrt{\rho})^2)\\
&+\left(2(\lambda+2\mu)\partial_3\log\sqrt{\rho}+2c_P\partial_3c_P\rho\right)c_P^{-1}J^{-1/2}c_P^{-1/2}(-\frac{1}{2}\rho^{-1/2}\partial_3\log\sqrt{\rho})\\
&+c_P^{-1/2}\sqrt{\rho}J^{-1/2}\partial_3\log\sqrt{\rho}\partial_3c_P+2c_P^{1/2}\sqrt{\rho}\partial_3 J^{-1/2}\partial_3\log\sqrt{\rho}\\
&+c_P^{-2}(\partial_3c_P)(\lambda+2\mu)J^{-1/2}c_P^{-1/2}\rho^{-1/2}\partial_3\log\sqrt{\rho}\\
&+2g^{\alpha\beta}c_S^2\rho\partial_\beta\log\sqrt{\rho}\left(\partial_\alpha J^{-1/2}c_P^{1/2}\rho^{-1/2}+\frac{1}{2}J^{-1/2}c_P^{-1/2}(\partial_\alpha c_P)\rho^{-1/2}-J^{-1/2}c_P^{1/2}\rho^{-1/2}\partial_\alpha\log\sqrt{\rho}\right)\\
&-2c_s(\partial_\beta c_S)\rho^{1/2}J^{-1/2}g^{\alpha\beta}(\partial_\alpha\log\sqrt{\rho})\\
&+\mu g^{\alpha\beta}\Big(-2\partial_\alpha J^{-1/2}c_P^{1/2}\rho^{-1/2}\partial_\beta\log\sqrt{\rho}+J^{-1/2}\rho^{-1/2}c_P^{-1/2}\partial_\alpha c_P\partial_\beta
\log\sqrt{\rho}\\
&\quad\quad+\rho^{-1/2}J^{-1/2}c_P^{1/2}(-\partial_\alpha\partial_\beta\log\sqrt{\rho}+(\partial_\alpha\log\sqrt{\rho})(\partial_\beta\log\sqrt{\rho}))\Big)\\
&+(c_P^2-2c_S^2)\rho\partial_3\log\sqrt{\rho} c_P^{-1}g^{\alpha\beta}\frac{\partial g^{\alpha\beta}}{\partial\tau}c_P^{-1/2}J^{-1/2}\rho^{-1/2}\\
&-3g^{\alpha\beta}\mu \partial_\beta c_PJ^{-1/2}c_P^{-1/2}\partial_\alpha\log\sqrt{\rho}\\
&-\frac{1}{2}g^{\alpha\beta}\frac{\partial g_{\alpha\beta}}{\partial\tau}(\lambda+2\mu) c_P^{-1}J^{-1/2}c_P^{-1/2}\rho^{-1/2}\partial_3\log\sqrt{\rho}+g^{\varepsilon\delta}\Gamma^\alpha_{\varepsilon\delta}\mu c_PJ^{-1/2}c_P^{-1/2}\rho^{-1/2}\partial_\alpha\log\sqrt{\rho}\\
&+\sqrt{\rho}\times\text{terms that do not depend on $\rho$}.
\end{align*}
The above formula simplifies to 
\begin{align*}
H_2=&-J^{-1/2}c_P^{1/2}\rho^{1/2}\partial^2_{33}\log\sqrt{\rho}-c_S^{2}\rho^{1/2}J^{-1/2}c_P^{1/2}g^{\alpha\beta}\partial_\alpha\partial_\beta\log\sqrt{\rho}\\
&-c_S^{2}\rho^{1/2}J^{-1/2}c_P^{1/2}g^{\alpha\beta}\partial_\alpha\log\sqrt{\rho}\partial_\beta\log\sqrt{\rho}-c_Pc_P^{-1/2}\rho^{1/2}J^{-1/2}(\partial_3\log\sqrt{\rho})^2\\
&+g^{\varepsilon\delta}\Gamma^\alpha_{\varepsilon\delta}c_S^2 c_P^{1/2}J^{-1/2}\rho^{1/2}\partial_\alpha\log\sqrt{\rho}+c_P^{-1/2}c_P^{-2}(3c_P^2)\rho^{1/2}J^{-1/2}\partial_3c_P\partial_3\log\sqrt{\rho}\\
&-3c_S^2c_P^{-1/2}\rho J^{-1/2}g^{\alpha\beta}\partial_\alpha c_P\partial_\beta\log\sqrt{\rho}-2c_s(\partial_\beta c_S)\rho^{1/2}J^{-1/2}g^{\alpha\beta}(\partial_\alpha\log\sqrt{\rho})\\
&+\frac{1}{2}g^{\alpha\beta}\frac{\partial g_{\alpha\beta}}{\partial\tau}(c_P^2-4c_S^2) c_P^{-1}J^{-1/2}c_P^{-1/2}\rho^{1/2}\partial_3\log\sqrt{\rho}\\
&+\sqrt{\rho}\times\text{terms that do not depend on $\rho$}.
\end{align*}
Then one can write
\begin{align*}
-\mathrm{i}H_1-H_2=&-\sqrt{\rho}(c_P^2-4c_S^2)c_P^{1/2}J^{-1/2}g^{\alpha\beta}\partial_\alpha\partial_\beta\log\sqrt{\rho}\\
&+\left(c_P^2-2c_S^2-\frac{2c_S^2(c_P^2-3c_S^2)}{c_P^2-c_S^2}\right)\sqrt{\rho}c_P^{1/2}J^{-1/2}g^{\alpha\beta}\partial_\alpha\log\sqrt{\rho}\partial_\beta\log\sqrt{\rho}\\
&+g^{\alpha\beta}\Gamma^\varepsilon_{\alpha\beta}(c_P^2-4c_S^2)c_P^{1/2}J^{-1/2}\partial_\varepsilon\log\sqrt{\rho}+J^{-1/2}c_P^{1/2}\rho^{1/2}\partial^2_{33}\log\sqrt{\rho}\\
&-\frac{1}{2}g^{\alpha\beta}\frac{\partial g_{\alpha\beta}}{\partial\tau}(c_P^2-4c_S^2) \rho^{1/2}c_P^{-1}J^{-1/2}c_P^{-1/2}\partial_3\log\sqrt{\rho}\\
&+c_P^{-1/2}(1+c_P^{-2}c_S^2)\rho^{1/2}J^{-1/2}\partial_3c_P\partial_3\log\sqrt{\rho}\\
&+\sqrt{\rho}\left(\frac{8c_S^3}{c_P^2-c_S^2}\right)c_P^{1/2}J^{-1/2}g^{\alpha\beta}\partial_\alpha\log\sqrt{\rho}\partial_\beta c_S-c_P^2c_P^{-1/2}\rho^{1/2}J^{-1/2}g^{\alpha\beta}\partial_\alpha\log\sqrt{\rho}\partial_\beta c_P\\
&+\sqrt{\rho}\times\text{terms that do not depend on $\rho$}.
\end{align*}

Using the formulas
\begin{equation}\label{fact1}
\nabla_{3}\nabla_{3}\log\sqrt{\rho}=\partial^2_{33}\log\sqrt{\rho}-c_P^{-1}\partial_3c_P\partial_3\log\sqrt{\rho}+c_Pg^{\alpha\beta}\partial_\beta c_P\partial_\alpha\log\sqrt{\rho},
\end{equation}
\begin{equation}\label{fact2}
\nabla_{\alpha}\nabla_{\beta}\log\sqrt{\rho}=\partial_\alpha\partial_\beta\log\sqrt{\rho}+\frac{1}{2}c_P^{-2}\frac{g_{\alpha\beta}}{\partial\tau}\partial_3\log\sqrt{\rho}-\Gamma^\varepsilon_{\alpha\beta}\partial_\varepsilon\log\sqrt{\rho},
\end{equation}
we obtain
\begin{align*}
&-\mathrm{i}H_1-H_2\\
=&\rho^{1/2}J^{-1/2}c_P^{1/2}\Big( -(c_P^2-4c_S^2)g^{\alpha\beta}\partial_\alpha\partial_\beta\log\sqrt{\rho}+(c_P^2-4c_S^2)g^{\alpha\beta}\Gamma^\varepsilon_{\alpha\beta}\partial_\varepsilon\log\sqrt{\rho}\\
&-(c_P^2-4c_S^2)\frac{1}{2}c_P^{-2}\frac{\partial g_{\alpha\beta}}{\partial\tau}g^{\alpha\beta}\partial_3\log\sqrt{\rho}-(c_P^2-4c_S^2)c_P^{-2}\partial^2_{33}\log\sqrt{\rho}+(c_P^2-4c_S^2)c_P^{-2}\partial_3c_P\partial_3\log\sqrt{\rho}\\
&-(c_P^2-4c_S^2)c_P^{-1}g^{\alpha\beta}\partial_\beta c_P\partial_\alpha\log\sqrt{\rho}\Bigg)\\
&+(2c_P^2-4c_S^2)\rho^{1/2}J^{-1/2}c_P^{1/2}\left(c_P^{-2}\partial^2_{33}\log\sqrt{\rho}-c_P^{-3}\partial_3 c_P\partial_3\log\sqrt{\rho}+c_P^{-1}g^{\alpha\beta}\partial_\beta c_P\partial_\alpha\log\sqrt{\rho}\right)\\
&+c_Pc_P^{-1/2}\rho^{1/2}J^{-1/2}(\partial_3\log\sqrt{\rho})^2\\
&+\sqrt{\rho}c_P^{1/2}J^{-1/2}\left(c_P^2-4c_S^2+\frac{4c_S^4}{c_P^2-c_S^2}\right)g^{\alpha\beta}\partial_\alpha\log\sqrt{\rho}\partial_\beta\log\sqrt{\rho}\\
&+\sqrt{\rho}c_P^{-3/2}J^{-1/2}\partial_3c_P^2\partial_3\log\sqrt{\rho}\\
&-c_P^{1/2}\rho^{1/2}J^{-1/2}g^{\alpha\beta}\partial_\alpha\log\sqrt{\rho}\partial_\beta c_P^2+\sqrt{\rho}\left(\frac{16c_S^3}{c_P^2-c_S^2}\right)c_P^{1/2}J^{-1/2}g^{\alpha\beta}\partial_\alpha\log\sqrt{\rho}\partial_\beta c_S\\
&+\sqrt{\rho}\times\text{terms that do not depend on $\rho$}.\\
=&\rho^{1/2}J^{-1/2}c_P^{1/2}\left(-(c_P^2-4c_S^2)\Delta\log\sqrt{\rho}\right)+(2c_P^2-4c_S^2)\rho^{1/2}J^{-1/2}c_P^{1/2}c_P^{-2}\nabla^2_{33}\log\sqrt{\rho}\\
&+\rho^{1/2}J^{-1/2}c_P^{1/2}\left(c_P^2-4c_S^2+\frac{4c_S^4}{c_P^2-c_S^2}\right)|\nabla\log\sqrt{\rho}|^2+4\rho^{1/2}J^{-1/2}c_P^{1/2}\frac{c_S^2}{c_P^2}\frac{c_P^2-2c_S^2}{c_P^2-c_S^2}|\partial_3\log\sqrt{\rho}|^2\\
&+\rho^{1/2}J^{-1/2}c_P^{1/2}\nabla\log\sqrt{\rho}\cdot\left(-\nabla c_P^2+\frac{8c_S^2}{c_P^2-c_S^2}\nabla c_S^2\right)\\
&+\rho^{1/2}J^{-1/2}c_P^{-3/2}\partial_3\log\sqrt{\rho}\left(2\partial_3 c_P^2-\frac{8c_S^2}{c_P^2-c_S^2}\partial_3c_S^2\right)\\
&+\sqrt{\rho}\times\text{terms that do not depend on $\rho$}.
\end{align*}
Inserting the relation $a_{13}=c_P^{1/2}\rho^{-1/2}J^{-1/2}B_P$ into the equation \eqref{transport_eq_amplitidue1P}, we prove the lemma.

\end{proof}

\begin{proof}[Proof of Lemma \ref{proofAS}]
With the choice of $\varphi=\varphi_S$, we calculate
\begin{align*}
&\partial_\alpha(\lambda a_{0k}\varphi_{;\ell}g^{k\ell})+\lambda a_{0k}\varphi_{;\ell}g^{k\ell}\partial_mg_{\alpha j} g^{jm}+\partial_m(\mu a_{0\alpha}\varphi_{;j}+\mu a_{0j}\varphi_{;\alpha})g^{jm}\\
=&\partial_3(\mu a_{0\alpha}\varphi_{;3})g^{33}\\
=&c_S^{-2}\partial_3(\mu a_{0\alpha}),\\
~\\
&\lambda a_{0k;\ell}g^{k\ell}\varphi_{;\alpha}+\mu(a_{0\alpha;j}+a_{0j;\alpha})\varphi_{;m}g^{jm}\\
=&c_S^{-2}\mu(\partial_3a_{0\alpha}-2\Gamma^\beta_{3\alpha}a_{0\beta})\\
=&c_S^{-2}\mu\partial_3a_{0\alpha}-c_S^{-2}\mu g^{\beta\varepsilon}\frac{\partial g_{\alpha\varepsilon}}{\partial\tau}a_{0\beta},\\
~\\
&\Gamma^n_{\alpha m}(\lambda a_{0k}\varphi_{;\ell}g^{k\ell}g_{nj}+\mu a_{0n}\varphi_{;j}+\mu a_{0j}\varphi_{;n})g^{jm}\\
=&g^{33}\Gamma^\varepsilon_{\alpha 3}\mu a_{0\varepsilon}+g^{\beta\varepsilon}\Gamma^3_{\alpha\beta}\mu a_{0\varepsilon}=0,\\
~\\
&\Gamma^n_{jm}(\lambda a_{0k}\varphi_{;\ell}g^{k\ell}g_{n\alpha}+\mu a_{0n}\varphi_{;\alpha}+\mu a_{0\alpha}\varphi_{;n})g^{jm}\\
=&(g^{33}\Gamma^3_{33}+g^{\beta\varepsilon}\Gamma^3_{\beta\varepsilon})\mu a_{0\alpha}\\
=&c_S^{-3}\partial_3c_S\mu a_{0\alpha}-\frac{1}{2}g^{\beta\varepsilon}\frac{\partial g_{\beta\varepsilon}}{\partial\tau}c_S^{-2}\mu a_{0\alpha}.
\end{align*}

As for the $P$-waves, we can also take $\mathbf{a}_j$ idependent of $t$. Using equation \eqref{second_eq} with $i=\alpha$, we obtain
\[
c_S^{-2}\left(2\mu\partial_3 a_{0\alpha}+(\partial_3\mu)a_{0\alpha}-\mu g^{\beta\varepsilon}\frac{\partial g_{\alpha\varepsilon}}{\partial\tau}a_{0\beta}+\frac{1}{2}g^{\beta\varepsilon}\frac{\partial g_{\beta\varepsilon}}{\partial\tau}c_S^{-2}\mu a_{0\alpha}-c_S^{-1}\partial_3c_S\mu a_{0\alpha}\right)=0.
\]
Introducing the vector $\eta$ such that $\eta_3=0$ and
\begin{equation}\label{equationeta}
\frac{\partial \eta_\alpha}{\partial \tau}-\frac{1}{2}g^{\beta\varepsilon}\frac{\partial g_{\alpha\beta}}{\partial\tau}\eta_\varepsilon+c_S^{-1}\partial_3c_S\eta_\alpha=0.
\end{equation}
and denote
\[
a_\alpha =A_Sc_S\eta_\alpha.
\]
If we consider the equation \eqref{equationeta} as an ODE along a geodesic $\gamma$ with respect to the metric $g_S$, we can see that $\eta$ is a parallel cotangent vector field along $\gamma$. So one can take $|\eta|^2_{g_S}=c_S^2g^{ij}\eta_i\eta_j=c_S^2g^{\alpha\beta}\eta_\alpha\eta_\beta=1$.
We obtain the following equation for $A_S$,
\[
c_S^{-2}\left(2\mu c_S\partial_3A_S+(\partial_3\mu)c_SA_S-\partial_3c_S\mu A_S+\frac{1}{2}\mu g^{\beta\varepsilon}\frac{\partial g_{\beta\varepsilon}}{\partial\tau}c_SA_S\right)=0.
\]
Similar to the proof of Lemma \ref{proofAP}, by solving the above ODE, we prove the lemma.
\end{proof}

\begin{proof}[Proof of Lemma \ref{proofBS}]
First we determine $a_{13}$. For this, we calculate
\begin{align*}
&\partial_3(\lambda a_{0k}\varphi_{;\ell}g^{k\ell})+\lambda a_{0k}\varphi_{;\ell}g^{k\ell}\partial_mg_{3 j} g^{jm}+\partial_m(\mu a_{03}\varphi_{;j}+\mu a_{0j}\varphi_{;3})g^{jm}\\
=&g^{\alpha\beta}\partial_\beta(\mu a_{0\alpha}),\\
~\\
&\lambda a_{0k;\ell}g^{k\ell}\varphi_{;3}+\mu(a_{03;j}+a_{0j;3})\varphi_{;m}g^{jm}\\
=&\lambda(\partial_\beta a_{0\alpha}-\Gamma^\varepsilon_{\alpha\beta}a_{0\varepsilon})g^{\alpha\beta}-\lambda\Gamma^\alpha_{33}a_{0\alpha}g^{33}+2g^{33}\mu(-\Gamma_{33}^\alpha a_{0\alpha})\\
=&\lambda g^{\alpha\beta}\partial_\beta a_{0\alpha}-g^{\alpha\beta}\lambda\Gamma^\varepsilon_{\alpha\beta}a_{0\varepsilon}+(\lambda+2\mu)c_S^{-1}g^{\alpha\beta}\partial_\beta c_Sa_{0\alpha},\\
~\\
&\Gamma^n_{3 m}(\lambda a_{0k}\varphi_{;\ell}g^{k\ell}g_{nj}+\mu a_{0n}\varphi_{;j}+\mu a_{0j}\varphi_{;n})g^{jm}\\
=&g^{33}\Gamma^\alpha_{33}\mu a_{0\alpha}+g^{\alpha\beta}\Gamma^3_{3\beta}\mu a_{0\alpha}\\
=&-c_S^{-1}g^{\alpha\beta}\partial_\beta c_S\mu a_{0\alpha}+g^{\alpha\beta}\mu c_S^{-1}\partial_\beta c_Sa_{0\alpha}=0,\\
~\\
&\Gamma^n_{jm}(\lambda a_{0k}\varphi_{;\ell}g^{k\ell}g_{n3}+\mu a_{0n}\varphi_{;3}+\mu a_{03}\varphi_{;n})g^{jm}\\
=&g^{jm}\Gamma_{jm}^\alpha\mu a_{0\alpha}\\
=&-c_S^{-1}\mu g^{\alpha\beta}\partial_\beta c_S a_{0\alpha}+g^{\beta\varepsilon}\Gamma^\alpha_{\beta\varepsilon}\mu a_{0\alpha}.
\end{align*}

Using equation \eqref{second_eq} with $i=3$, we obtain
\[
\begin{split}
\rho a_{13}-(\lambda+2\mu)c_S^{-2}a_{13}\\
+\mathrm{i}\Big(g^{\alpha\beta}\partial_\beta(\mu a_{0\alpha})+\lambda g^{\alpha\beta}\partial_\beta a_{0\alpha}-\lambda g^{\alpha\beta}\Gamma^\varepsilon_{\alpha\beta}a_{0\varepsilon}+(\lambda+2\mu)c_S^{-1}g^{\alpha\beta}\partial_\beta c_Sa_{0\alpha}\\
+c_S^{-1}\mu g^{\alpha\beta}\partial_\beta c_S a_{0\alpha}-g^{\beta\varepsilon}\Gamma^\alpha_{\beta\varepsilon}\mu a_{0\alpha}\Big)=0,
\end{split}
\]
which simplifies to
\[
\begin{split}
\rho a_{13}-(\lambda+2\mu)c_S^{-2}a_{13}\\
+\mathrm{i}\Big((\lambda+\mu) g^{\alpha\beta}\partial_\beta a_{0\alpha}+g^{\alpha\beta}(\partial_\beta\mu) a_{0\alpha}+(\lambda+3\mu)c_S^{-1}g^{\alpha\beta}\partial_\beta c_Sa_{0\alpha}-(\lambda+\mu)g^{\beta\varepsilon}\Gamma^\alpha_{\beta\varepsilon} a_{0\alpha}\Big)=0.
\end{split}
\]
So
\[
\begin{split}
a_{13}=&\frac{\mathrm{i}}{c_P^2c_S^{-2}-1}\left(\frac{\lambda+\mu}{\rho} g^{\alpha\beta}\partial_\beta a_{0\alpha}+g^{\alpha\beta}\frac{\partial_\beta\mu}{\rho} a_{0\alpha}+\frac{\lambda+3\mu}{\rho}c_S^{-1}g^{\alpha\beta}\partial_\beta c_Sa_{0\alpha}-\frac{\lambda+\mu}{\rho}g^{\beta\varepsilon}\Gamma^\alpha_{\beta\varepsilon} a_{0\alpha}\right)\\
=&\frac{\mathrm{i}}{c_P^2c_S^{-2}-1}\Bigg((c_P^2-c_S^2)g^{\beta\varepsilon}\partial_\beta a_{0\varepsilon}+g^{\alpha\beta}\partial_\beta c_S^2 a_{0\varepsilon}+2g^{\beta\varepsilon}c_S^2\partial_\beta\log\sqrt{\rho}a_{0\varepsilon}+(c_P^2+c_S^2)c_S^{-1}g^{\beta\varepsilon}\partial_\beta c_S a_{0\varepsilon}\\
&\quad\quad\quad\quad-(c_P^2-c_S^2)g^{\beta\varepsilon}\Gamma^\alpha_{\beta\varepsilon} a_{0\alpha}\Bigg).
\end{split}
\]
To continue, we calculate
\begin{align*}
\partial_\alpha a_{13}=&\frac{\mathrm{i}}{c_P^2c_S^{-2}-1}\Bigg((c_P^2-c_S^2)g^{\beta\varepsilon}\partial_\alpha\partial_\beta a_{0\varepsilon}+\partial_\alpha(c_P^2-c_S^2)g^{\beta\varepsilon}\partial_\beta a_{0\varepsilon}+(c_P^2-c_S^2)\partial_\alpha g^{\beta\varepsilon}\partial_\beta a_{0\varepsilon}\\
&\quad \quad +g^{\beta\varepsilon}(\partial_\beta c_S^2)\partial_\alpha a_{0\varepsilon}+2g^{\beta\varepsilon}c_S^2\partial_\beta\log\sqrt{\rho}\partial_\alpha a_{0\varepsilon}+2\partial_\alpha g^{\beta\varepsilon}c_S^2\partial_\beta\log\sqrt{\rho} a_{0\varepsilon}\\
&\quad \quad+2 g^{\beta\varepsilon}\partial_\alpha c_S^2\partial_\beta\log\sqrt{\rho} a_{0\varepsilon}+2g^{\beta\varepsilon}c_S^2\partial_\alpha \partial_\beta\log\sqrt{\rho} a_{0\varepsilon}\\
&\quad\quad +(c_P^2+c_S^2)c_S^{-1}g^{\beta\varepsilon}\partial_\beta c_S\partial_\alpha a_{0\varepsilon}-(c_P^2+c_S^2)g^{\beta\varepsilon}\Gamma^\delta_{\beta\varepsilon}\partial_\alpha a_{0\delta}\\
&\quad\quad -\frac{1}{c_S^2}\frac{c_S^2\partial_\alpha c_P^2+c_P^2\partial_\alpha c_S^2}{c_P^2-c_S^2}(c_P^2-c_S^2)g^{\beta\varepsilon}\partial_\beta a_{0\varepsilon}-\frac{2}{c_S^2}\frac{c_S^2\partial_\alpha c_P^2+c_P^2\partial_\alpha c_S^2}{c_P^2-c_S^2}g^{\beta\varepsilon}c_S^2\partial_\beta\log\sqrt{\rho} a_{0\varepsilon}\Bigg)\\
&+\frac{1}{\sqrt{\rho}}\times\text{terms that do not depend on $\rho$}\\
=&\frac{\mathrm{i}}{c_P^2c_S^{-2}-1}\Bigg((c_P^2-c_S^2)g^{\beta\varepsilon}\partial_\alpha\partial_\beta a_{0\varepsilon}-\frac{c_P^2+c_S^2}{c_S^2}\partial_\alpha(c_S^2)g^{\beta\varepsilon}\partial_\beta a_{0\varepsilon}+(c_P^2-c_S^2)\partial_\alpha g^{\beta\varepsilon}\partial_\beta a_{0\varepsilon}\\
&\quad \quad +2g^{\beta\varepsilon}c_S^2\partial_\beta\log\sqrt{\rho}\partial_\alpha a_{0\varepsilon}+2\partial_\alpha g^{\beta\varepsilon}c_S^2\partial_\beta\log\sqrt{\rho} a_{0\varepsilon}\\
&\quad \quad+2g^{\beta\varepsilon}c_S^2\partial_\alpha \partial_\beta\log\sqrt{\rho} a_{0\varepsilon}\\
&\quad\quad +(c_P^2+3c_S^2)c_S^{-1}g^{\beta\varepsilon}\partial_\beta c_S\partial_\alpha a_{0\varepsilon}-(c_P^2-c_S^2)g^{\beta\varepsilon}\Gamma^\delta_{\beta\varepsilon}\partial_\alpha a_{0\delta}\\
&\quad\quad +2 g^{\beta\varepsilon}\partial_\alpha c_S^2\partial_\beta\log\sqrt{\rho} a_{0\varepsilon}-2\frac{c_S^2\partial_\alpha c_P^2+c_P^2\partial_\alpha c_S^2}{c_P^2-c_S^2}g^{\beta\varepsilon}\partial_\beta\log\sqrt{\rho} a_{0\varepsilon}\Bigg)\\
&+\frac{1}{\sqrt{\rho}}\times\text{terms that do not depend on $\rho$}.
\end{align*}
We proceed with the tedious calculation
\begin{align*}
&\partial_\alpha(\lambda a_{1k}\varphi_{;\ell}g^{k\ell})+\lambda a_{1k}\varphi_{;\ell}g^{k\ell}\partial_mg_{\alpha j} g^{jm}+\partial_m(\mu a_{1\alpha}\varphi_{;j}+\mu a_{1j}\varphi_{;\alpha})g^{jm}\\
=&\partial_\alpha(\lambda a_{13}g^{33})+\lambda a_{13}g^{33}g^{\beta\varepsilon}\partial_\beta g_{\alpha\varepsilon}+\partial_3(\mu a_{1\alpha})g^{33}\\
=&\partial_\alpha(\lambda a_{13}c_S^{-2})+\lambda a_{13}c_S^{-2}g^{\beta\varepsilon}\partial_\beta g_{\alpha\varepsilon}+\partial_3(\mu a_{1\alpha})c_S^{-2},\\
~\\
&\lambda a_{1k;\ell}g^{k\ell}\varphi_{;\alpha}+\mu(a_{1\alpha;j}+a_{1j;\alpha})\varphi_{;m}g^{jm}\\
=&\mu(\partial_3 a_{1\alpha})g^{33}+\mu(\partial_\alpha a_{13})g^{33}-2\mu\Gamma^i_{\alpha 3}a_{1i}g^{33}\\
=&\mu c_S^{-2}\partial_3a_{1\alpha}+\mu c_S^{-2}\partial_\alpha a_{13}-2\mu c_S^{-2}(\frac{1}{2}g^{\beta\varepsilon}\frac{\partial g_{\alpha\varepsilon}}{\partial\tau}a_{1\beta}+c_S^{-1}\partial_\alpha c_S a_{13})\\
=&\mu c_S^{-2}\partial_3a_{1\alpha}-\mu c_S^{-2}g^{\beta\varepsilon}\frac{\partial g_{\alpha\varepsilon}}{\partial\tau}a_{1\beta}+\mu c_S^{-2}\partial_\alpha a_{13}-2\mu c_S^{-3}\partial_\alpha c_S a_{13},\\
~\\
&\Gamma^n_{\alpha m}(\lambda a_{1k}\varphi_{;\ell}g^{k\ell}g_{nj}+\mu a_{1n}\varphi_{;j}+\mu a_{1j}\varphi_{;n})g^{jm}\\
=&\Gamma_{\alpha n}^n\lambda a_{13}g^{33}+\Gamma^n_{\alpha 3}\mu a_{1n}g^{33}+\Gamma^3_{\alpha m}\mu a_{1j}g^{jm}\\
=&(\Gamma_{\alpha 3}^3+\Gamma_{\alpha \beta}^\beta )\lambda a_{13}g^{33}+\Gamma^3_{\alpha 3}\mu a_{13}g^{33}+\Gamma^\beta_{\alpha 3}\mu a_{1\beta}g^{33}+\Gamma^3_{\alpha 3}\mu a_{13}g^{33}+\Gamma^3_{\alpha \beta}\mu a_{1\gamma}g^{\beta\gamma}\\
=&(c_S^{-1}\partial_\alpha c_S+\Gamma^\beta_{\alpha\beta})\lambda c_S^{-2}a_{13}+2\mu c_S^{-3}\partial_\alpha c_S a_{13},\\
~\\
&\Gamma^n_{jm}(\lambda a_{1k}\varphi_{;\ell}g^{k\ell}g_{n\alpha}+\mu a_{1n}\varphi_{;\alpha}+\mu a_{1\alpha}\varphi_{;n})g^{jm}\\
=&g^{jm}\Gamma^\beta_{jm}\lambda a_{13}g^{33}g_{\alpha\beta}+(g^{33}\Gamma^3_{33}+g^{\beta\varepsilon}\Gamma^3_{\beta\varepsilon})\mu a_{1\alpha}\\
=&-\lambda c_S^{-3}\partial_\alpha c_S a_{13}+c_S^{-2}g^{\varepsilon\delta}g_{\alpha\beta}\varepsilon^\beta_{\varepsilon\delta}\lambda a_{13}+c_S^{-3}\partial_3c_S\mu a_{1\alpha}-\frac{1}{2}g^{\beta\varepsilon}\frac{\partial g_{\beta\varepsilon}}{\partial\tau}c_S^{-2}\mu a_{1\alpha},
\end{align*}
and
\begin{align*}
&\partial_m(\lambda a_{0k;\ell}g^{k\ell}g_{\alpha j}+\mu a_{0\alpha ;j}+\mu a_{0j;\alpha})g^{jm}\\
=&g^{\beta\varepsilon}\partial_\beta(\lambda g^{33}(-\Gamma_{33}^\delta a_{0\delta})g_{\alpha \varepsilon}+\lambda(\partial_\varepsilon a_{0\delta}-\Gamma_{\kappa\delta}^\eta a_{0\eta} )g^{\delta\kappa}g_{\alpha \varepsilon})\\
&+g^{jm}\partial_m(\mu(\partial_ja_{0\alpha}-\Gamma^\beta_{\alpha j}a_{0\beta}))+g^{jm}\partial_m(\mu(\partial_\alpha a_{0j}-\Gamma^\beta_{\alpha j}a_{0\beta}))\\
=&g^{\beta\varepsilon}\partial_\beta(\lambda(\partial_\kappa a_{0\delta})g^{\delta\kappa}g_{\alpha\varepsilon})-g^{\beta\varepsilon}\partial_\beta(\lambda g^{33}\Gamma^\delta_{33}a_{0\delta}g_{\alpha\gamma}+\lambda\Gamma^\eta_{\kappa\delta}a_{0\eta}g^{\delta\kappa}g_{\alpha\varepsilon})\\
&+g^{33}\partial_3(\mu\partial_3a_{0\alpha})+g^{\beta\varepsilon}\partial_\beta(\mu\partial_\varepsilon a_{0\alpha})+g^{\beta\varepsilon}\partial_\varepsilon(\mu\partial_\alpha a_{0\beta})-2g^{33}\partial_3(\mu\Gamma^\beta_{\alpha 3}a_{0\beta})-2g^{\delta\varepsilon}\partial_\delta(\mu \Gamma^\beta_{\alpha\varepsilon} a_{0\beta})\\
=&\partial_\alpha \lambda g^{\delta\varepsilon}\partial_\varepsilon a_{0\delta}+\lambda g^{\delta\varepsilon}\partial_{\alpha}\partial_\varepsilon a_{0\delta}+\lambda \partial_\varepsilon a_{0\delta}\partial_\alpha g^{\delta\varepsilon}+\lambda g^{\beta\kappa}g^{\delta\varepsilon}\partial_\beta g_{\alpha\kappa}\partial_\varepsilon a_{0\delta}\\
&+(\partial_\alpha\lambda)c_S^{-1}g^{\delta\beta}\partial_\beta c_Sa_{0\delta}+\lambda c_S^{-1}g^{\delta\beta}\partial_\beta c_S\partial_\alpha a_{0\delta}-\partial_\alpha\lambda g^{\delta\varepsilon}\Gamma^\eta_{\varepsilon\delta}a_{0\eta}-\lambda g^{\delta\varepsilon}\Gamma^\eta_{\varepsilon\delta}\partial_\alpha a_{0\eta}\\
&+c_S^{-2}\partial_3\mu\partial_3a_{0\alpha}+c_S^{-2}\mu\partial^2_{33}a_{0\alpha}+g^{\beta\varepsilon}\partial_\beta\mu\partial_\varepsilon a_{0\alpha}+\mu g^{\beta\varepsilon}\partial_\beta\partial_\varepsilon a_{0\alpha}+g^{\beta\varepsilon}\partial_\varepsilon\mu\partial_\alpha a_{0\beta}+\mu g^{\beta\varepsilon}\partial_\alpha\partial_\varepsilon a_{0\beta}\\
&-2c_S^{-2}\mu \Gamma^\beta_{\alpha 3} \partial_3a_{0\beta}-2g^{\delta\varepsilon}\mu \Gamma^\beta_{\alpha\varepsilon} \partial_\delta a_{0\beta}-2c_S^{-2}\partial_3\mu\Gamma^\beta_{\alpha 3}a_{0\beta}-2g^{\delta\varepsilon}\partial_\delta\mu\Gamma^\beta_{\alpha\varepsilon}a_{0\beta},\\
~\\
&\Gamma_{\alpha m}^n(\lambda a_{0k;\ell}g^{k\ell}g_{nj}+\mu a_{0n;j}+\mu a_{0j;n})g^{jm}\\
=&\Gamma^n_{\alpha n}\lambda g^{\beta\gamma}\partial_\gamma a_{0\beta}+g^{jm}\Gamma^\beta_{\alpha m}\mu\partial_j a_{0\beta}+g^{\beta\gamma}\mu\Gamma^n_{\alpha \gamma}\partial_n a_{0\beta}+\sqrt{\rho}\times\text{terms that do not depend on $\rho$}\\
=&(\Gamma^3_{3\alpha}+\Gamma^\delta_{\alpha\delta})\lambda g^{\beta\varepsilon}\partial_\varepsilon a_{0\beta}+g^{33}\Gamma^\beta_{\alpha 3}\mu\partial_3 a_{0\beta}+g^{\varepsilon\delta}\Gamma^\beta_{\alpha\varepsilon}\mu\partial_\delta a_{0\beta}+g^{\beta\varepsilon}\mu\Gamma^3_{\alpha\varepsilon}\partial_3 a_{0\beta}+g^{\beta\varepsilon}\Gamma^\delta_{\alpha\varepsilon}\mu\partial_\delta a_{0\beta}\\
&+\sqrt{\rho}\times\text{terms that do not depend on $\rho$}\\
=&(\Gamma^3_{3\alpha}+\Gamma^\delta_{\alpha\delta})\lambda g^{\beta\varepsilon}\partial_\varepsilon a_{0\beta}+g^{\varepsilon\delta}\Gamma^\beta_{\alpha\varepsilon}\mu\partial_\delta a_{0\beta}+g^{\beta\varepsilon}\Gamma^\delta_{\alpha\varepsilon}\mu\partial_\delta a_{0\beta}+\sqrt{\rho}\times\text{terms that do not depend on $\rho$},\\
~\\
&\Gamma_{jm}^n(\lambda a_{0k;\ell}g^{k\ell}g_{n\alpha}+\mu a_{0n;\alpha}+\mu a_{0\alpha;n})g^{jm}\\
=&g^{jm}\Gamma^\beta_{jm}g_{\alpha\beta}\lambda \partial_\delta a_{0\varepsilon}g^{\varepsilon\delta}+g^{jm}\Gamma^\beta_{jm}\mu\partial_\alpha a_{0\beta}+g^{jm}\Gamma^\beta_{jm}\mu\partial_\beta a_{0\alpha}+g^{jm}\Gamma^3_{jm}\mu\partial_3 a_{0\alpha}\\
&+\sqrt{\rho}\times\text{terms that do not depend on $\rho$}\\
=&{g^{33}\Gamma^\beta_{33}g_{\alpha\beta}\lambda \partial_\delta a_{0\varepsilon}g^{\varepsilon\delta}}+{g^{\eta\varepsilon}\Gamma^\beta_{\eta\varepsilon}g_{\alpha\beta}\lambda \partial_\delta a_{0\varepsilon}g^{\varepsilon\delta}}+\mu(g^{33}\Gamma^\beta_{33}+ g^{\eta\varepsilon}\Gamma_{\eta\varepsilon}^\beta)(\partial_\alpha a_{0\beta}+\partial_\beta a_{0\alpha})\\
&+(g^{33}\Gamma^3_{33}+g^{\beta\varepsilon}\Gamma^3_{\beta\varepsilon})\mu\partial_3 a_{0\alpha}+\sqrt{\rho}\times\text{terms that do not depend on $\rho$}.
\end{align*}

Using equation \eqref{third_eq} with $i=\alpha$, we obtain
\begin{equation}\label{eq_swave_amplitude1}
\mathrm{i}c_S^{-2}\left(2\mu\partial_3 a_{1\alpha}+(\partial_3\mu)a_{1\alpha}-\mu g^{\beta\varepsilon}\frac{\partial g_{\alpha\varepsilon}}{\partial\tau}a_{1\beta}+\frac{1}{2}g^{\beta\varepsilon}\frac{\partial g_{\beta\varepsilon}}{\partial\tau}c_S^{-2}\mu a_{1\alpha}-c_S^{-1}\partial_3c_S\mu a_{1\alpha}\right)=\rho H_\alpha,
\end{equation}
where
\begin{align*}
H_\alpha=&{(c_P^2-c_S^2)g^{\beta\varepsilon}\partial_\alpha\partial_\beta a_{0\varepsilon}}-\frac{c_P^2+c_S^2}{c_S^2}\partial_\alpha(c_S^2)g^{\beta\varepsilon}\partial_\beta a_{0\varepsilon}+(c_P^2-c_S^2)\partial_\alpha g^{\beta\varepsilon}\partial_\beta a_{0\varepsilon}\\
& +2g^{\beta\varepsilon}c_S^2\partial_\beta\log\sqrt{\rho}\partial_\alpha a_{0\varepsilon}+2\partial_\alpha g^{\beta\varepsilon}c_S^2\partial_\beta\log\sqrt{\rho} a_{0\varepsilon}\\
&+2g^{\beta\varepsilon}c_S^2\partial_\alpha \partial_\beta\log\sqrt{\rho} a_{0\varepsilon}\\
& +(c_P^2+3c_S^2)c_S^{-1}g^{\beta\varepsilon}\partial_\beta c_S\partial_\alpha a_{0\varepsilon}-{(c_P^2-c_S^2)g^{\beta\varepsilon}\Gamma^\delta_{\beta\varepsilon}\partial_\alpha a_{0\delta}}\\
&+2 g^{\beta\varepsilon}\partial_\alpha c_S^2\partial_\beta\log\sqrt{\rho} a_{0\varepsilon}-2\frac{c_S^2\partial_\alpha c_P^2+c_P^2\partial_\alpha c_S^2}{c_P^2-c_S^2}g^{\beta\varepsilon}\partial_\beta\log\sqrt{\rho} a_{0\varepsilon}\\
&-2c_P^2c_S^{-1}\partial_\alpha c_S\frac{1}{c_P^2-c_S^2}\left((c_P^2-c_S^2)g^{\beta\varepsilon}\partial_\beta a_{0\varepsilon}+2g^{\beta\varepsilon}c_S^2\partial_\beta\log\sqrt{\rho}a_{0\varepsilon}\right)\\
&+\partial_\alpha(c_P^2-2c_S^2)\frac{1}{c_P^2-c_S^2}\left({(c_P^2-c_S^2)g^{\beta\varepsilon}\partial_\beta a_{0\varepsilon}}+2g^{\beta\varepsilon}c_S^2\partial_\beta\log\sqrt{\rho}a_{0\varepsilon}\right)\\
&+\partial_\alpha\log\sqrt{\rho}\frac{2(c_P^2-2c_S^2)}{c_P^2-c_S^2}\left({\frac{\lambda+\mu}{\rho} g^{\beta\varepsilon}\partial_\beta a_{0\varepsilon}}+g^{\varepsilon\beta}\frac{\partial_\beta\mu}{\rho} a_{0\varepsilon}+\frac{\lambda+3\mu}{\rho}c_S^{-1}g^{\beta\varepsilon}\partial_\beta c_Sa_{0\varepsilon}-{\frac{\lambda+\mu}{\rho}g^{\beta\varepsilon}\Gamma^\delta_{\beta\varepsilon} a_{0\delta}}\right)\\
&-{\rho^{-1}\partial_\alpha \lambda g^{\delta\varepsilon}\partial_\varepsilon a_{0\delta}}-{(c_P^2-2c_S^2) g^{\delta\varepsilon}\partial_{\alpha}\partial_\varepsilon a_{0\delta}}-(c_P^2-2c_S^2) \partial_\varepsilon a_{0\delta}\partial_\alpha g^{\delta\varepsilon}\\
&-\rho^{-1}\partial_\alpha\lambda c_S^{-1}g^{\delta\beta}\partial_\beta c_Sa_{0\delta}-(c_P^2-2c_S^2) c_S^{-1}g^{\delta\beta}\partial_\beta c_S\partial_\alpha a_{0\delta}+{\rho^{-1}\partial_\alpha\lambda g^{\delta\varepsilon}\Gamma^\eta_{\varepsilon\delta}a_{0\eta}}+{(c_P^2-2c_S^2) g^{\delta\varepsilon}\Gamma^\eta_{\varepsilon\delta}\partial_\alpha a_{0\eta}}\\
&-c_S^{-2}\rho^{-1}\partial_3\mu\partial_3a_{0\alpha}-\partial^2_{33}a_{0\alpha}-\rho^{-1}g^{\beta\varepsilon}\partial_\beta\mu\partial_\varepsilon a_{0\alpha}-c_S^2 g^{\beta\varepsilon}\partial_\beta\partial_\varepsilon a_{0\alpha}-g^{\beta\varepsilon}\rho^{-1}\partial_\varepsilon\mu\partial_\alpha a_{0\beta}-{\mu g^{\beta\varepsilon}\partial_\alpha\partial_\varepsilon a_{0\beta}}\\
&+2 \Gamma^\beta_{\alpha 3} \partial_3a_{0\beta}+2g^{\delta\varepsilon}c_S^2 \Gamma^\beta_{\alpha\varepsilon} \partial_\delta a_{0\beta}+\partial_3\rho g^{\beta\varepsilon}\frac{\partial g_{\alpha\varepsilon}}{\partial \tau}a_{0\beta}+2c_S^2g^{\delta\varepsilon}\partial_\delta\rho\Gamma^\beta_{\alpha\varepsilon}a_{0\beta}\\
&+c_S^{-1}\partial_3 c_S\partial_3 a_{0\alpha}-\frac{1}{2}g^{\beta\varepsilon}\frac{\partial g_{\beta\varepsilon}}{\partial\tau}\partial_3 a_{0\alpha}\\
&+g^{\varepsilon\delta}\Gamma^\beta_{\alpha\varepsilon}c_S^2\partial_\delta a_{0\beta}+g^{\beta\varepsilon}\Gamma^\delta_{\alpha\varepsilon}c_S^2\partial_\delta a_{0\beta}\\
&+c_S^2(-c_S^{-1}g^{\beta\varepsilon}\partial_\varepsilon c_S)(\partial_\alpha a_{0\beta}+\partial_\beta a_{0\alpha})+ c_S^2 g^{\eta\varepsilon}\Gamma_{\eta\varepsilon}^\beta({\partial_\alpha a_{0\beta}}+\partial_\beta a_{0\alpha})\\
&+\frac{1}{\sqrt{\rho}}\times\text{terms that do not depend on $\rho$}.
\end{align*}
Here we have also used the identity \eqref{Christofel_identity}.\\

Tedious calculations give
\begin{align*}
H_\alpha=&-\partial^2_{33}a_{0\alpha}-c_S^{-2}\rho^{-1}\partial_3\mu\partial_3a_{0\alpha}+2 \Gamma^\beta_{\alpha 3} \partial_3a_{0\beta}\\
&-g^{\beta\varepsilon}\rho^{-1}\partial_\beta\mu\partial_\varepsilon a_{0\alpha}-c_S^2 g^{\beta\varepsilon}\partial_\beta\partial_\varepsilon a_{0\alpha}-\rho^{-1}g^{\beta\varepsilon}\partial_\varepsilon\mu\partial_\alpha a_{0\beta}\\
&+2g^{\beta\varepsilon}c_S^2\partial_\beta\log\sqrt{\rho}\partial_\alpha a_{0\varepsilon}+2g^{\beta\varepsilon}c_S^2\partial_\alpha \partial_\beta\log\sqrt{\rho} a_{0\varepsilon}\\
&+2(c_P^2-2c_S^2)\partial_\alpha\log\sqrt{\rho}\frac{1}{c_P^2-c_S^2}g^{\varepsilon\beta}\frac{\partial_\beta\mu}{\rho} a_{0\varepsilon}\\
&-\frac{2c_P^2+c_S^2}{c_S^2}\partial_\alpha(c_S^2)g^{\beta\varepsilon}\partial_\beta a_{0\varepsilon}-\partial_\alpha (c_P^2-2c_S^2)g^{\delta\varepsilon }\partial_\varepsilon a_{0\delta}\\
&+(c_P^2-c_S^2)\partial_\alpha g^{\beta\varepsilon}\partial_\beta a_{0\varepsilon}+2\partial_\alpha g^{\beta\varepsilon}c_S^2\partial_\beta\log\sqrt{\rho} a_{0\varepsilon}-(c_P^2-2c_S^2)\partial_\varepsilon a_{0\delta}\partial_\alpha g^{\delta\varepsilon}\\
&+3g^{\delta\varepsilon}c_S^2 \Gamma^\beta_{\alpha\varepsilon} \partial_\delta a_{0\beta}+g^{\beta\varepsilon}\Gamma^\delta_{\alpha\varepsilon}c_S^2\partial_\delta a_{0\beta}\\
&+2\partial_3\log\sqrt{\rho}g^{\beta\varepsilon}\frac{\partial g_{\alpha\varepsilon}}{\partial\tau}a_{0\beta}+4c_S^2g^{\delta\varepsilon}\partial_\delta\log\sqrt{\rho}\Gamma^\beta_{\alpha\varepsilon}a_{0\beta}\\
&+(-2\partial_\alpha c_S^2-\frac{8c_S^2\partial c_S^2}{c_P^2-c_S^2})g^{\beta\varepsilon}\partial_\beta\log\sqrt{\rho} a_{0\varepsilon}\\
&+2(c_P^2-2c_S^2)\partial_\alpha\log\sqrt{\rho}\frac{1}{c_P^2-c_S^2}\frac{\lambda+3\mu}{\rho}c_S^{-1}g^{\beta\varepsilon}\partial_\beta c_Sa_{0\varepsilon}-\rho^{-1}\partial_\alpha\lambda c_S^{-1}g^{\delta\beta}\partial_\beta c_Sa_{0\delta}\\
&-c_S\partial_\varepsilon c_{S}g^{\beta\varepsilon}\partial_\beta a_{0\alpha}+\mu g^{\eta\varepsilon}\Gamma_{\eta\varepsilon}^\beta\partial_\beta a_{0\alpha}\\
&+c_S^{-3}\partial_3 c_S\mu\partial_3 a_{0\alpha}-\frac{1}{2}g^{\beta\varepsilon}\frac{\partial g_{\beta\varepsilon}}{\partial\tau}\partial_3 a_{0\alpha}\\
&+4c_Sg^{\beta\varepsilon}\partial_\beta c_S\partial_\alpha a_{0\varepsilon}\\
&+\frac{1}{\sqrt{\rho}}\times\text{terms that do not depend on $\rho$}.
\end{align*}

Recalling $a_{0\alpha} =A_Sc_S\eta_\alpha$, we have

\begin{align*}
H_\alpha=&-\partial^2_{33}A_S c_S\eta_\alpha{-2\partial_3A_S\partial_3 (c_S\eta_\alpha})-c_S^{-1}\partial_3c_S^2\partial_3A_S\eta_\alpha-2\partial_3\log\sqrt{\rho}\partial_3A_Sc_S\eta_\alpha\\
&{-2\partial_3(c_S\log\sqrt{\rho}A_S\partial_3\eta_\alpha})+g^{\beta\varepsilon}\frac{\partial g_{\alpha\varepsilon}}{\partial \tau}\partial_3A_Sc_S\eta_\beta\\
&-g^{\beta\varepsilon}\partial_\beta c_S^2\partial_\varepsilon A_Sc_S\eta_\alpha-2g^{\beta\varepsilon}c_S^3\partial_\beta\log\sqrt{\rho}\partial_\varepsilon A_S\eta_\alpha-{2g^{\beta\varepsilon}c_S^2\partial_\beta\log\sqrt{\rho} A_S\partial_\varepsilon(c_S\eta_\alpha})\\
&-\mu g^{\beta\varepsilon}\partial_\beta\partial_\varepsilon A_Sc_S\eta_\alpha-{2\mu g^{\beta\varepsilon}\partial_\beta A_S\partial_\varepsilon(c_S\eta_\alpha})\\
&-g^{\beta\varepsilon}\partial_\varepsilon c_S^2\partial_\alpha A_Sc_S\eta_\beta-{2g^{\beta\varepsilon}c_S^2\partial_\varepsilon\log\sqrt{\rho}\partial_\alpha A_Sc_S\eta_\beta}-{2g^{\beta\varepsilon}c_S^2\partial_\varepsilon\log\sqrt{\rho} A_S\partial_\alpha(c_S\eta_\beta})\\
&+{2g^{\beta\varepsilon}c_S^3\partial_\beta\log\sqrt{\rho}\partial_\alpha A_S \eta_\varepsilon}+{2g^{\beta\varepsilon}c_S^2\partial_\beta\log\sqrt{\rho} A_S \partial_\alpha(c_S\eta_\varepsilon})+2g^{\beta\varepsilon}c_S^3\partial_\alpha \partial_\beta\log\sqrt{\rho} a_{0\varepsilon}\\
&+2\partial_\alpha\log\sqrt{\rho}\frac{c_P^2-2c_S^2}{c_P^2-c_S^2}g^{\varepsilon\beta}\partial_\beta c_S^2 A_Sc_S\eta_{\varepsilon}+4\partial_\alpha\log\sqrt{\rho}\frac{c_P^2-2c_S^2}{c_P^2-c_S^2}g^{\varepsilon\beta} c_S^2\partial_\beta\log\sqrt{\rho} A_Sc_S\eta_{\varepsilon}\\
&-\frac{2c_P^2+c_S^2}{c_S^2}\partial_\alpha(c_S^2)g^{\beta\varepsilon}\partial_\beta A_Sc_S \eta_{\varepsilon}-\partial_\alpha (c_P^2-2c_S^2)g^{\delta\varepsilon }\partial_\varepsilon A_Sc_S\eta_\delta\\
&+(c_P^2-c_S^2)\partial_\alpha g^{\beta\varepsilon}\partial_\beta A_S c_S\eta_{\varepsilon}+2\partial_\alpha g^{\beta\varepsilon}c_S^2\partial_\beta\log\sqrt{\rho}c_S A_S\eta_{\varepsilon}\\
&-\lambda \partial_\varepsilon A_S c_S\eta_{\delta}\partial_\alpha g^{\delta\varepsilon}+3g^{\delta\varepsilon}\mu \Gamma^\beta_{\alpha\varepsilon} \partial_\delta A_Sc_S \eta_{\beta}+g^{\beta\varepsilon}\Gamma^\delta_{\alpha\varepsilon}\mu\partial_\delta A_Sc_S \eta_{\beta}\\
&+(-2\partial_\alpha c_S^2-\frac{8c_S^2\partial c_S^2}{c_P^2-c_S^2})g^{\beta\varepsilon}\partial_\beta\log\sqrt{\rho} A_Sc_S \eta_{\varepsilon}\\
&+2(c_P^2-2c_S^2)\partial_\alpha\log\sqrt{\rho}\frac{1}{c_P^2-c_S^2}\frac{\lambda+3\mu}{\rho}c_S^{-1}g^{\beta\varepsilon}\partial_\beta c_S^2 A_S\eta_{\varepsilon}-2\partial_\alpha\log\sqrt{\rho}c_Sg^{\delta\beta}\partial_\beta c^2_S A_S\eta_{\delta}\\
&-c_S\partial_\varepsilon c_{S}g^{\beta\varepsilon}\partial_\beta A_Sc_S \eta_{\alpha}+c_S^2 g^{\eta\varepsilon}\Gamma_{\eta\varepsilon}^\beta\partial_\beta A_Sc_S \eta_{\alpha}\\
&+c_S^{-3}\partial_3 c_S\mu c_S\partial_3 A_S \eta_\alpha-\frac{1}{2}g^{\beta\varepsilon}c_S^{-2}\frac{\partial g_{\beta\varepsilon}}{\partial\tau}\mu c_S\partial_3A_S \eta_{\alpha}\\
&+2\partial_3\log\sqrt{\rho}g^{\beta\varepsilon}\frac{\partial g_{\alpha\varepsilon}}{\partial\tau}c_SA_S\eta_{\beta}+4c_S^2g^{\delta\varepsilon}\partial_\delta\log\sqrt{\rho}\Gamma^\beta_{\alpha\varepsilon}c_SA_S\eta_{\beta}\\
&+4c_Sg^{\beta\varepsilon}\partial_\beta c_S\partial_\alpha c_SA_S \eta_{\varepsilon}\\
&+\frac{1}{\sqrt{\rho}}\times\text{terms that do not depend on $\rho$}.
\end{align*}
Upon using
\[
\partial A_S=-J^{-1/2}c_S^{-1/2}\rho^{-1/2}\partial\log\sqrt{\rho}-\frac{1}{2}c_S^{-3/2}\partial_3 c_S J^{-1/2}\rho^{-1/2}-\frac{1}{2}J^{-3/2}\partial Jc_S^{-1/2}\rho^{-1/2},
\]
\[
\begin{split}
\partial^2_{33} A_S=&-J^{-1/2}c_S^{-1/2}\rho^{-1/2}\partial^2_{3}\log\sqrt{\rho}+J^{-1/2}c_S^{-1/2}\rho^{-1/2}(\partial_{3}\log\sqrt{\rho})^2\\
&+c_S^{-3/2}\partial_3 c_S J^{-1/2}\rho^{-1/2}\partial_3\log\sqrt{\rho}+J^{-3/2}\partial_3 Jc_S^{-1/2}\rho^{-1/2}\partial_3\log\sqrt{\rho}\\
&+\frac{1}{\sqrt{\rho}}\times\text{terms that do not depend on $\rho$},
\end{split}
\]
and
\[
\begin{split}
g^{\beta\varepsilon}\partial_\beta\partial_\varepsilon A_S=&J^{-1/2}c_S^{-1/2}g^{\beta\varepsilon}\partial_\beta\partial_\varepsilon \rho^{-1/2}+2g^{\beta\varepsilon}\partial_\beta J^{-1/2}c_S^{-1/2}\partial_\varepsilon\rho^{-1/2}+2J^{-1/2}g^{\beta\varepsilon}\partial_\beta c_S^{-1/2}\partial_\varepsilon\rho^{-1/2}\\
=&-J^{-1/2}c_S^{-1/2}\rho^{-1/2}g^{\beta\varepsilon}\partial_\alpha\partial_\beta\log\sqrt{\rho}+J^{-1/2}c_S^{-1/2}\rho^{-1/2}g^{\beta\varepsilon}(\partial_\alpha\log\sqrt{\rho})(\partial_\beta\log\sqrt{\rho})\\
&+g^{\beta\varepsilon}J^{-3/2}\partial_\beta J c_P^{-1/2}\rho^{-1/2}\partial_\varepsilon\log\sqrt{\rho}+J^{-1/2}\rho^{-1/2}g^{\beta\varepsilon}c_S^{-3/2}\partial_\beta c_S\partial_\varepsilon\log\sqrt{\rho}\\
&+\frac{1}{\sqrt{\rho}}\times\text{terms that do not depend on $\rho$},
\end{split}
\]
we have
\[
\rho H=\rho^{1/2}c_S^{1/2}J^{-1/2}H',
\]
where
\begin{align*}
H'_\alpha=&\partial^2_{3}\log\sqrt{\rho}\eta_\alpha-(\partial_{3}\log\sqrt{\rho})^2\eta_\alpha-c_S^{-1}\partial_3 c_S \partial_3\log\sqrt{\rho}\eta_\alpha-{J^{-1}\partial_3 J\partial_3\log\sqrt{\rho}\eta_\alpha}+c_S^{-2}\partial_3c_S^2\partial_3\log\sqrt{\rho}\eta_\alpha\\
&-2\partial_3\log\sqrt{\rho}\left(-\partial_3\log\sqrt{\rho}-\frac{1}{2}c_S^{-1}\partial_3 c_S -{\frac{1}{2}J^{-1}\partial_3 J}\right)\eta_\alpha+g^{\beta\varepsilon}\frac{\partial g_{\alpha\varepsilon}}{\partial \tau}\partial_3\log\sqrt{\rho}\eta_\beta\\
&+g^{\beta\varepsilon}\partial_\beta c_S^2\partial_\varepsilon \log\sqrt{\rho}\eta_\alpha-2g^{\beta\varepsilon}c_S^2\partial_\beta\log\sqrt{\rho}\left(-\partial_\varepsilon\log\sqrt{\rho}-\frac{1}{2}c_S^{-1}\partial_\varepsilon c_S-{\frac{1}{2}J^{-1}\partial_\varepsilon J}\right)\eta_\alpha\\
&+c_S^2g^{\beta\varepsilon}\partial_\varepsilon\partial_\beta\log\sqrt{\rho}\eta_\alpha-c_S^2g^{\beta\varepsilon}(\partial_\varepsilon\log\sqrt{\rho})(\partial_\beta\log\sqrt{\rho})\eta_\alpha-{c_S^2g^{\beta\varepsilon}J^{-1}\partial_\beta J \partial_\varepsilon\log\sqrt{\rho}\eta_\alpha}\\
&-c_S^2g^{\beta\varepsilon}c_S^{-1}\partial_\beta c_S\partial_\varepsilon\log\sqrt{\rho}\eta_\alpha\\
&+g^{\beta\varepsilon}\partial_\varepsilon c_S^2\partial_\alpha \log\sqrt{\rho}\eta_\beta+2g^{\beta\varepsilon}c_S^2\partial_\alpha \partial_\beta\log\sqrt{\rho}\eta_\varepsilon\\
&+2(c_P^2-2c_S^2)\partial_\alpha\log\sqrt{\rho}\frac{1}{c_P^2-c_S^2}g^{\varepsilon\beta}\partial_\beta c_S^2 \eta_{\varepsilon}+4(c_P^2-2c_S^2)\partial_\alpha\log\sqrt{\rho}\frac{1}{c_P^2-c_S^2}g^{\varepsilon\beta} c_S^2\partial_\beta\log\sqrt{\rho} \eta_{\varepsilon}\\
&+\frac{2c_P^2+c_S^2}{c_S^2}\partial_\alpha(c_S^2)g^{\beta\varepsilon}\partial_\beta\log\sqrt{\rho}  \eta_{\varepsilon}+\partial_\alpha(c_P^2- 2c_S^2)g^{\delta\varepsilon }\partial_\varepsilon\log\sqrt{\rho}\eta_\delta\\
&-c_S^2\partial_\alpha g^{\beta\varepsilon}\partial_\beta \log\sqrt{\rho} \eta_{\varepsilon}+2\partial_\alpha g^{\beta\varepsilon}c_S^2\partial_\beta\log\sqrt{\rho} \eta_{\varepsilon}+g^{\delta\varepsilon}c_S^2 \Gamma^\beta_{\alpha\varepsilon} \partial_\delta  \log\sqrt{\rho} \eta_{\beta}-g^{\beta\varepsilon}\Gamma^\delta_{\alpha\varepsilon}c_S^2\partial_\delta  \log\sqrt{\rho} \eta_{\beta}\\
&+(-2\partial_\alpha c_S^2-\frac{8c_S^2\partial_\alpha c_S^2}{c_P^2-c_S^2})g^{\beta\varepsilon}\partial_\beta\log\sqrt{\rho}  \eta_{\varepsilon}\\
&+2(c_P^2-2c_S^2)\partial_\alpha\log\sqrt{\rho}\frac{1}{c_P^2-c_S^2}(c_P^2+c_S^2)c_S^{-1}g^{\beta\varepsilon}\partial_\beta c_S \eta_{\varepsilon}-2\partial_\alpha\log\sqrt{\rho}c_Sg^{\delta\beta}\partial_\beta c_S \eta_{\delta}\\
&+c_S\partial_\varepsilon c_{S}g^{\beta\varepsilon}\partial_\beta \log\sqrt{\rho}  \eta_{\alpha}-c_S^2 g^{\eta\varepsilon}\Gamma_{\eta\varepsilon}^\beta\partial_\beta \log\sqrt{\rho}  \eta_{\alpha}\\
&-c_S^{-3}\partial_3 c_S c_S^2\partial_3 \log\sqrt{\rho}\eta_\alpha+\frac{1}{2}g^{\beta\varepsilon}c_S^{-2}\frac{\partial g_{\beta\varepsilon}}{\partial\tau}c_S^2\partial_3 \log\sqrt{\rho}\eta_\alpha-4c_S^2 g^{\beta\varepsilon}\partial_\beta c_S\partial_\alpha \log\sqrt{\rho}\eta_\varepsilon\\
&+\text{terms that do not depend on $\rho$}.
\end{align*}

Note that
\[
\begin{split}
&\partial_\alpha g^{\beta\varepsilon}c_S^2\partial_\beta\log\sqrt{\rho} \eta_{\varepsilon}+g^{\delta\varepsilon}c_S^2 \Gamma^\beta_{\alpha\varepsilon} \partial_\delta  \log\sqrt{\rho} \eta_{\beta}-g^{\beta\varepsilon}\Gamma^\delta_{\alpha\varepsilon}c_S^2\partial_\delta  \log\sqrt{\rho} \eta_{\beta}\\
=&c_S^2\left(-g^{\beta \delta}\partial_\alpha g_{\delta\varepsilon}g^{\varepsilon\kappa}+\frac{1}{2}g^{\beta\delta}g^{\kappa\varepsilon}(\partial_\alpha g_{\delta\varepsilon}+\partial_\delta g_{\alpha\varepsilon}-\partial_\varepsilon g_{\alpha\delta})-\frac{1}{2}g^{\kappa\delta}g^{\beta\varepsilon}(\partial_\alpha g_{\delta\varepsilon}+\partial_\delta g_{\alpha\varepsilon}-\partial_\varepsilon g_{\alpha\delta})\right)\partial_\beta\log\sqrt{\rho} \eta_{\kappa}\\
=&-c_S^2g^{\kappa\delta}g^{\beta\varepsilon}(\partial_\alpha g_{\delta\varepsilon}+\partial_\delta g_{\alpha\varepsilon}-\partial_\varepsilon g_{\alpha\delta})\partial_\beta\log\sqrt{\rho} \eta_{\kappa}\\
=&-2c_S^2g^{\varepsilon\delta}\Gamma^{\beta}_{\alpha\delta}\partial_\beta\log\sqrt{\rho} \eta_{\varepsilon}.
\end{split}
\]
Notice also that
\[
c_S^{2}|\eta|^2=c^2_{S}g^{\alpha\beta}\eta_\alpha\eta_\beta=c^2_{S}g_{\alpha\beta}\eta^\alpha\eta^\beta=1.
\]
Using also \eqref{fact1} and \eqref{fact2}, we then have
\begin{align*}
&\langle H',\eta\rangle_{g_S}=c_S^{2}g^{\alpha\beta}H'_\alpha\eta_{\beta}\\
=&\partial^2_{3}\log\sqrt{\rho}+(\partial_{3}\log\sqrt{\rho})^2+c_S^2g^{\beta\varepsilon}(\partial_\varepsilon\log\sqrt{\rho})(\partial_\beta\log\sqrt{\rho})+c_S^2g^{\beta\varepsilon}\partial_\varepsilon\partial_\beta\log\sqrt{\rho}+\\
&+2c_S^4\partial_\alpha \partial_\beta\log\sqrt{\rho}\eta^\alpha\eta^\beta+4(c_P^2-2c_S^2)\frac{c_2^4}{c_P^2-c_S^2}\partial_\alpha\log\sqrt{\rho}\partial_\beta\log\sqrt{\rho} \eta^{\beta} \eta^{\alpha}-c_S^2 g^{\eta\varepsilon}\Gamma_{\eta\varepsilon}^\beta\partial_\beta \log\sqrt{\rho} \\
&+\frac{1}{2}g^{\beta\varepsilon}\frac{\partial g_{\beta\varepsilon}}{\partial\tau}\partial_3 \log\sqrt{\rho}+c_S^{2}\frac{\partial g_{\alpha\varepsilon}}{\partial \tau}\partial_3\log\sqrt{\rho}\eta^\varepsilon\eta^\alpha-2c_S^4\Gamma^{\beta}_{\alpha\delta}\partial_\beta\log\sqrt{\rho} \eta^{\delta}\eta^\alpha\\
&+c_S^{-1}\partial_3c_S\partial_3\log\sqrt{\rho}+3c_Sg^{\beta\varepsilon}\partial_\beta\log\sqrt{\rho}\partial_\varepsilon c_S\\
&+3(c_P^2-c_S^2)\partial_\alpha c_S^2 \partial_\beta\log\sqrt{\rho}\eta^\alpha\eta^\beta+c_S^{2}\partial_\alpha c_P^2 \partial_\beta\log\sqrt{\rho}\eta^\alpha\eta^\beta-\frac{12c_S^4}{c_P^2-c_S^2}\partial_\alpha c_S^2 \partial_\beta\log\sqrt{\rho}\eta^\alpha\eta^\beta\\
&+\text{terms that do not depend on $\rho$}\\
=&\partial^2_{3}\log\sqrt{\rho}-c_S^{-1}\partial_3c_S\partial_3\log\sqrt{\rho}+c_Sg^{\beta\varepsilon}\partial_\beta\log\sqrt{\rho}\partial_\varepsilon c_S+c_S^2g^{\beta\varepsilon}\partial_\varepsilon\partial_\beta\log\sqrt{\rho}\\
&+\frac{1}{2}g^{\beta\varepsilon}\frac{\partial g_{\beta\varepsilon}}{\partial\tau}\partial_3 \log\sqrt{\rho}-c_S^2 g^{\eta\varepsilon}\Gamma_{\eta\varepsilon}^\beta\partial_\beta \log\sqrt{\rho}\\
&+2c_S^4\partial_\alpha \partial_\beta\log\sqrt{\rho}\eta^\beta\eta^\alpha+c_S^{2}\frac{\partial g_{\alpha\varepsilon}}{\partial \tau}\partial_3\log\sqrt{\rho}\eta^\varepsilon\eta^\alpha-c_S^42\Gamma^{\beta}_{\alpha\delta}\partial_\beta\log\sqrt{\rho} \eta^{\delta}\eta^\alpha\\
&+c_S^2g^{\beta\varepsilon}(\partial_\varepsilon\log\sqrt{\rho})(\partial_\beta\log\sqrt{\rho})+(\partial_{3}\log\sqrt{\rho})^2\\
&+4(c_P^2-2c_S^2)\frac{c_S^4}{c_P^2-c_S^2} \partial_\alpha\log\sqrt{\rho}\partial_\beta\log\sqrt{\rho} \eta^{\beta} \eta^{\alpha}\\
&+c_S^{2}\partial_3c_S^2\partial_3\log\sqrt{\rho}+g^{\beta\varepsilon}\partial_\beta\log\sqrt{\rho}\partial_\varepsilon c_S^2\\
&+3(c_P^2-c_S^2)\partial_\alpha c_S^2 \partial_\beta\log\sqrt{\rho}\eta^\alpha\eta^\beta+c_S^{2}\partial_\alpha c_P^2 \partial_\beta\log\sqrt{\rho}\eta^\alpha\eta^\beta-\frac{12c_S^4}{c_P^2-c_S^2}\partial_\alpha c_S^2 \partial_\beta\log\sqrt{\rho}\eta^\alpha\eta^\beta\\
&+\text{terms that do not depend on $\rho$}\\
=&c_S^2\Delta\log\sqrt{\rho}+2c_S^4\nabla_\alpha\nabla_\beta\log\sqrt{\rho}\eta^\alpha\eta^\beta+c_S^2|\nabla\log\sqrt{\rho}|^2+4c_S^4\frac{c_P^2-2c_S^2}{c_P^2-c_S^2}\nabla_\alpha\log\sqrt{\rho}\nabla_\beta\log\sqrt{\rho}\eta^\alpha\eta^\beta\\
&+\nabla\log\sqrt{\rho}\cdot\nabla c_S^2+c_S^4\nabla_\beta\log\sqrt{\rho}\left(2\frac{c_P^2-c_S^2}{c_S^4}\nabla_\alpha c_S^2+c_S^{-2}\nabla_\alpha c_P^2-\frac{12}{c_P^2-c_S^2}\nabla_\alpha c_S^2\right)\eta^\alpha\eta^\beta\\
&+\text{terms that do not depend on $\rho$}\\
=&:c^4_Sg^{\alpha\varepsilon}g^{\beta\delta}N_{\alpha\beta}\eta_\varepsilon\eta_\delta+C\\
=&g_S^{\alpha\varepsilon}g_S^{\beta\delta}N_{\alpha\beta}\eta_\varepsilon\eta_\delta+C,
\end{align*}
where $C$ does not depend on $\rho$. Recall that $g_S=c_S^{-2}g_E$, and consequently
\[
g_S^{\alpha\beta}=c_S^2g^{\alpha\beta}.
\]
One should be careful in the above calculation that $\eta^\alpha\partial_{x^\alpha}$ is tangent vector corresponding to the cotangent vector $\eta$ w.r.t. the Euclidean metric $g_E$. The corresponding tangent vector w.r.t. the metric $g_S$ should be $c_S^2\eta^\alpha\partial_{x^\alpha}$. Note that the second order tensor $N$ can be written in a coordinate-free way
\[
\begin{split}
N=&c_S^2\Delta\log\sqrt{\rho}g_S+2\nabla^2\log\sqrt{\rho}+c_S^2|\nabla\log\sqrt{\rho}|^2g_S+4\frac{c_P^2-2c_S^2}{c_P^2-c_S^2}\nabla\log\sqrt{\rho}\otimes\nabla\log\sqrt{\rho}\\
&+(\nabla\log\sqrt{\rho}\cdot\nabla c_S^2) g_S+\nabla\log\sqrt{\rho}\otimes\left(2\frac{c_P^2-c_S^2}{c_S^4}\nabla c_S^2+c_S^{-2}\nabla c_P^2-\frac{12}{c_P^2-c_S^2}\nabla c_S^2\right).
\end{split}
\]

Multiplying both sides of \eqref{eq_swave_amplitude1} by $c_S^{2}\eta^\alpha$, we get
\[
\begin{split}
&\mathrm{i}c_S^{-2}\left(2\mu\partial_3 a_{1\alpha}+(\partial_3\mu)a_{1\alpha}-\mu g^{\beta\varepsilon}\frac{\partial g_{\alpha\varepsilon}}{\partial\tau}a_{1\beta}+\frac{1}{2}g^{\beta\varepsilon}\frac{\partial g_{\beta\varepsilon}}{\partial\tau}c_S^{-2}\mu a_{1\alpha}-c_S^{-1}\partial_3c_S\mu a_{1\alpha}\right)c_S^2\eta^\alpha \\
=&(\langle N,\eta\otimes\eta \rangle_{g_S} +C)J^{-1/2}c_S^{1/2}\rho^{1/2}.
\end{split}
\]
 Denote $D_\tau$ to be the covariant derivative along the geodesic. We consider the above equation along $\gamma$.
 If we denote
 \[
 a_{1i}=J^{-1/2}c_S^{1/2}\rho^{-1/2}\xi_i,
 \]
 with a vector field $\xi=(\xi_1,\xi_2,\xi_3)$ along $\gamma$,
the above equation can be written as
 \[
D_\tau \langle \xi,\eta\rangle_{g_S}= \langle D_\tau\xi,\eta\rangle_{g_S}=-\frac{\mathrm{i}}{2}\langle N,\eta\otimes\eta\rangle_{g_S}+C
 \]
 along $\gamma$. The proof is complete if one notices that $\langle \xi,\eta\rangle_{g_S}=B_S$.
\end{proof}

{\bf Acknowledgement.}
The author thanks the anonymous referees for very useful comments that helped improve this paper significantly.  \\

\noindent \textbf{Data Availability} This article has no associated data.\\

\noindent \textbf{Conflict of interest} The authors have no conflict of interest to declare that are relevant to the content of this article.
\bibliographystyle{abbrv}
\bibliography{biblio}

\end{document}